\documentclass[notitlepage]{article}
\usepackage{babel} 															
\usepackage{titling}
\usepackage{xcolor} 														
\usepackage[top=3cm,bottom=3cm,left=3cm,right=3cm]{geometry}
\usepackage{fancyhdr}
\usepackage{mathtools}
\usepackage{unicode-math}
\usepackage{amsthm}
\usepackage[Rejne]{fncychap}
\usepackage{titletoc}
\usepackage{minitoc}
\usepackage{titlesec}
\usepackage{enumitem}
\usepackage{imakeidx}
\usepackage[colorlinks=true,linkcolor=blue,citecolor=orange,urlcolor=purple,hyperindex]{hyperref}
\usepackage{cleveref}
\usepackage{array}
\usepackage[framemethod=tikz]{mdframed}
\usepackage{xpatch}
\usepackage[all]{xy}

\newcommand{\batonA}{\msansA}
\newcommand{\batonD}{\msansD}
\newcommand{\batonE}{\msansE}
\newcommand{\batonP}{\msansP}
\newcommand{\batonQ}{\msansQ}

\newcommand{\cA}{{\mathcal{A}}}
\newcommand{\cD}{{\mathcal{D}}}
\newcommand{\cF}{{\mathcal{F}}}
\newcommand{\cG}{{\mathcal{G}}}
\newcommand{\cH}{{\mathcal{H}}}
\newcommand{\cI}{{\mathcal{I}}}
\newcommand{\cJ}{{\mathcal{J}}}
\newcommand{\cK}{{\mathcal{K}}}
\newcommand{\cL}{{\mathcal{L}}}
\newcommand{\cM}{{\mathcal{M}}}
\newcommand{\cN}{{\mathcal{N}}}
\newcommand{\cO}{{\mathcal{O}}}
\newcommand{\cP}{{\mathcal{P}}}
\newcommand{\cQ}{{\mathcal{Q}}}

\newcommand{\fC}{{\mathfrak{C}}}
\newcommand{\fD}{{\mathfrak{D}}}
\newcommand{\fE}{{\mathfrak{E}}}
\newcommand{\mfA}{{\mathfrak{a}}}
\newcommand{\mfC}{{\mathfrak{c}}}
\newcommand{\mfL}{{\mathfrak{l}}}
\newcommand{\mfM}{{\mathfrak{m}}}

\newcommand{\eA}{{\mathbb{A}}}
\newcommand{\eD}{{\mathbb{D}}}
\newcommand{\eN}{{\mathbb{N}}}
\newcommand{\eQ}{{\mathbb{Q}}}
\newcommand{\eR}{{\mathbb{R}}}

\newcommand{\et}	{\text{ and }}

\definecolor{couleurdeff}	{rgb}	{0.09, 0.45, 0.27}
\definecolor{couleurtheo}	{rgb}	{0.77, 0.12, 0.23}
\definecolor{couleurlem}	{rgb}	{0.99, 0.76, 0.0} 
\definecolor{couleurcor}	{rgb}	{0.54, 0.17, 0.89} 
\definecolor{couleurprop}	{rgb}	{1.0, 0.44, 0.37}

\titleformat{\section}[frame]{\normalfont}{\filright\footnotesize\enspace Section \thesection\enspace}{8pt}{\Large\bfseries\filcenter}
\titleformat{\subsection}{\large\bfseries}{\thesubsection.}{0.5em}{}
  
\newlength{\itemsepa}
\newlength{\topsepa}
\newlength{\labelwidtha}
\newlength{\leftmargina}
\setitemize{itemsep=\itemsepa ,topsep=\topsepa, label = --}
\setenumerate{itemsep=\itemsepa ,topsep=\topsepa, leftmargin=\leftmargina, labelwidth=\labelwidtha}
\setenumerate[2]{itemsep=\itemsepa ,topsep=\topsepa, label = \alph*)}

\newcommand{\sautlistebis}{\leavevmode\vspace{-0.0cm}}

\makeatletter
\newcommand{\skipitems}[1]{%
  \addtocounter{\@enumctr}{#1}%
}
\makeatother

\makeatletter
\newtheoremstyle{definition}
	{2pt}
	{2pt}
	{\slshape}
	{0pt}
	{}
	{\leavevmode}
	{\newline}
	{\textcolor{couleurdeff}{\thmname{\textbf{#1}}}\;\thmnumber{\textbf{#2}} \thmnote{--- \textbf{[#3]}}}
\theoremstyle{definition}

\newtheorem{deff_internal}{Définition}[subsection]
\newenvironment{deff}[1][]%
	{\begin{mdframed}[
		innerlinecolor=couleurdeff,
		outerlinecolor=couleurdeff,
		innertopmargin=-0.1\baselineskip,
		innerbottommargin=0.0\baselineskip,
		innerleftmargin=5pt,
		innerrightmargin=1pt,
		outerlinewidth=0.7pt,
		innerlinewidth=0.7pt,
		middlelinewidth=0.7pt,
		leftmargin=10pt,
		hidealllines=true,
		leftline=true,
		middlelinecolor=white,
		bottomline=false,
		footnoteinside=true,
		skipbelow=10\baselineskip,
	]\begin{deff_internal}[#1]}%
	{\end{deff_internal}\end{mdframed}}

\newtheoremstyle{lemme}
	{2pt}
	{2pt}
	{\slshape}
	{0pt}
	{}
	{}
	{0pt}
	{\textcolor{couleurlem}{\thmname{\textbf{#1}}}\;\thmnumber{\textbf{#2}} \thmnote{--- \textbf{[#3]}\newline}}
\theoremstyle{lemme}

\newtheorem{lem_internal}[deff_internal]{Lemme}
\newenvironment{lem}[1][]%
	{\begin{mdframed}[
		innerlinecolor=couleurlem,
		outerlinecolor=couleurlem,
		innertopmargin=-0.1\baselineskip,
		innerbottommargin=0.0\baselineskip,
		innerleftmargin=5pt,
		innerrightmargin=1pt,
		outerlinewidth=0.7pt,
		innerlinewidth=0.7pt,
		middlelinewidth=0.7pt,
		leftmargin=10pt,
		hidealllines=true,
		leftline=true,
		middlelinecolor=white,
		bottomline=false,
		footnoteinside=true,
		skipbelow=10\baselineskip,
	]\begin{lem_internal}[#1]}%
	{\end{lem_internal}\end{mdframed}}

\newtheoremstyle{proposition}
	{2pt}
	{2pt}
	{\slshape}
	{0pt}
	{}
	{}
	{0pt}
	{\textcolor{couleurprop}{\thmname{\textbf{#1}}}\;\thmnumber{\textbf{#2}} \thmnote{--- \textbf{[#3]}\newline}}
\theoremstyle{proposition}

\newtheorem{prop_internal}[deff_internal]{Proposition}
\newenvironment{prop}[1][]%
	{\begin{mdframed}[
		innerlinecolor=couleurprop,
		outerlinecolor=couleurprop,
		innertopmargin=-0.1\baselineskip,
		innerbottommargin=0.0\baselineskip,
		innerleftmargin=5pt,
		innerrightmargin=1pt,
		outerlinewidth=0.7pt,
		innerlinewidth=0.7pt,
		middlelinewidth=0.7pt,
		leftmargin=10pt,
		hidealllines=true,
		leftline=true,
		middlelinecolor=white,
		bottomline=false,
		footnoteinside=true,
		skipbelow=10\baselineskip,
	]\begin{prop_internal}[#1]}%
	{\end{prop_internal}\end{mdframed}}

\newtheoremstyle{theoreme}
	{2pt}
	{2pt}
	{\slshape}
	{0pt}
	{}
	{}
	{0pt}
	{\textcolor{couleurtheo}{\thmname{\textbf{#1}}}\;\thmnumber{\textbf{#2}} \thmnote{--- \textbf{[#3]}\newline}}
\theoremstyle{theoreme}

\newtheorem{theo_internal}[deff_internal]{Théorème}
\newenvironment{theo}[1][]%
	{\begin{mdframed}[
		innerlinecolor=couleurtheo,
		outerlinecolor=couleurtheo,
		innertopmargin=-0.1\baselineskip,
		innerbottommargin=0.0\baselineskip,
		innerleftmargin=5pt,
		innerrightmargin=1pt,
		outerlinewidth=0.7pt,
		innerlinewidth=0.7pt,
		middlelinewidth=0.7pt,
		leftmargin=10pt,
		hidealllines=true,
		leftline=true,
		middlelinecolor=white,
		bottomline=false,
		footnoteinside=true,
		skipbelow=10\baselineskip,
	]\begin{theo_internal}[#1]}%
	{\end{theo_internal}\end{mdframed}}

\newtheoremstyle{corollaire}
	{2pt}
	{2pt}
	{\slshape}
	{0pt}
	{}
	{}
	{0pt}
	{\textcolor{couleurcor}{\thmname{\textbf{#1}}}\;\thmnumber{\textbf{#2}} \thmnote{--- \textbf{[#3]}\newline}}
\theoremstyle{corollaire}

\newtheorem{cor_internal}[deff_internal]{Corollaire}
	{\begin{mdframed}[
		innerlinecolor=couleurcor,
		outerlinecolor=couleurcor,
		innertopmargin=-0.1\baselineskip,
		innerbottommargin=0.0\baselineskip,
		innerleftmargin=5pt,
		innerrightmargin=1pt,
		outerlinewidth=0.7pt,
		innerlinewidth=0.7pt,
		middlelinewidth=0.7pt,
		leftmargin=10pt,
		hidealllines=true,
		leftline=true,
		middlelinecolor=white,
		bottomline=false,
		footnoteinside=true,
		skipbelow=10\baselineskip,
	]\begin{cor_internal}[#1]}%
	{\end{cor_internal}\end{mdframed}}

\newtheoremstyle{remarques}
	{\topsep}
	{\topsep}
	{\slshape}
	{0pt}
	{}
	{}
	{0pt}
	{\thmname{\textbf{#1}}\;\thmnumber{\textbf{#2}} \thmnote{--- \textbf{[#3]}\newline}}
\theoremstyle{remarques}

\newtheorem{voc}[deff_internal]{Vocabulaire}
\newtheorem{att}[deff_internal]{Attention}
\newtheorem{nott}[deff_internal]{Notation}

\newtheoremstyle{demonstration}
	{0pt}
	{\topsep}
	{\normalfont\footnotesize}
	{0pt}
	{}
	{}
	{0pt}
	{\thmname{\textbf{#1}}\;\thmnumber{\textbf{#2}} \thmnote{--- \textbf{[#3]}\newline}}
\theoremstyle{demonstration}

\newtheorem{dem}	[deff_internal]{Démonstration}

\newtheoremstyle{approfondissement}
	{2pt}
	{2pt}
	{\small\normalfont%
	\addtolength{\@totalleftmargin}{4.3em} %
	\addtolength{\linewidth}{-4.3em} %
	\parshape 1 4.3em \linewidth} %
	{0pt}
	{}
	{}
	{0pt}
	{\thmname{\textbf{#1}}\;\thmnumber{\textbf{#2}} \thmnote{--- \textbf{[#3]}\newline}}
\theoremstyle{approfondissement}

\makeatother

\newcommand{\PM}					{\scriptscriptstyle}

\newcommand{\srl}			[1]		{\overline{#1}}

\newcommand{\BP}			[1]		{\mathopen{}\left({#1}\right)\mathclose{}}						
\newcommand{\BC}			[1]		{\mathopen{}\left\langle {#1} \right\rangle\mathclose{}}		

\def\cartesien{
\ar@{-}[]+L+<-4pt,-2pt>;[]+LD+<-4pt,-4pt>
\ar@{-}[]+D+<-2pt,-4pt>;[]+LD+<-4pt,-4pt>}
\def\cartesienbis{
\ar@{-}[]+R+<4pt,-2pt>;[]+RD+<4pt,-4pt>
\ar@{-}[]+D+<2pt,-4pt>;[]+RD+<4pt,-4pt>}

\makeatletter
\newcommand{\superimpose}[2]{{%
  \ooalign{%
    \hfil$\m@th#1\@firstoftwo#2$\hfil\cr
    \hfil$\m@th#1\@secondoftwo#2$\hfil\cr
  }%
}}
\makeatother
\newcommand{\superpose}[2]{\mathpalette\superimpose{{#1}{#2}}}

\newcommand{\lseg}					{\lBrack}
\newcommand{\rseg}					{\rBrack}

\newcommand{\no}					{n{}\textsuperscript{o}}

\newcommand{\absNP}			[2][]	{\left| #2 \right|{\vphantom{K}}_{\PM #1}^{\PM\!\!\times}} 		
\newcommand{\adh}			[2][]	{\srl{#2}{}^{\PM #1}} 												
\newcommand{\Abhy}			[2][]	{\batonA\BP{#2}} 													
\newcommand{\abs}			[2][]	{\left| #2 \right|_{\PM #1}} 										
\newcommand{\AnnRes}		[2][]	{
			\ifx o#1 				{#2}^{\circ}							\else
			\ifx v#1 				{#2}^{\circ\circ}						\else
			\ifx i#1 				{#2}^{\PM\sim}							\else
			\ifx r#1 				{#2}^{\PM\equiv}						\else
			\ifx t#1 				{#2}^{\PM\approx}						\else
			\mathrm{Illegal~option}									\fi\fi\fi\fi\fi} 					
\newcommand{\Algbr}			[2][]	{{#2}{}^{\mfA\mfL}}													

\newcommand{\BLup}			[3][]	{\batonE_{\PM #2}\BP{#3}}						
\newcommand{\BLupD}			[3][]	{\batonD_{\PM #2}\BP{#3}}						
\newcommand{\BorRelB}		[2][]	{\partial\BP{#2}}								
\newcommand{\BijL}			[2][]	{\textrm{Bij}\BP{#2}} 							

\newcommand{\Conc}			[2][]	{#2^{\mfC}}										
\newcommand{\cOf}			[2][]	{\cO^{\PM #1}_{#2}}								
\newcommand{\cOfg}			[3][]	{\cO^{\PM #1}_{#2,#3}}							
\newcommand{\CrpRes}		[2][]	{\kappa\BP{#2}}									
\newcommand{\CrpResC}		[2][]	{\cH\BP{#2}}									
\newcommand{\CrpResCR}		[2][]	{\tilde{\cH}\BP{#2}}							
\newcommand{\Coprofl}		[2][]	{\textrm{codepth}_{\PM #1}\BP{#2}} 			
\newcommand{\cosep}			[2][]	{\widehat{#2}}									
\newcommand{\CoSep}			[2][]	{\superpose{\vphantom{#2}}{\smash{\cosep{#2}}}} 	

\newcommand{\dimAn}			[3][]	{\textrm{Adim}_{\PM #2}^{\PM #1}\BP{#3}}		
\newcommand{\dimCtr}		[2][]	{\textrm{cdim}_{\PM #1}\BP{#2}}				
\newcommand{\dimKr}			[2][]	{\textrm{dim}_{#1}\BP{#2}}						
\newcommand{\DiagApplAn}	[3][\DEF]	{
			\ifx \DEF #1 	\fD_{\PM #2,#3}^{\PM\textbf{elem}} 		\else
			\fD_{\PM #2,#3}^{\PM\leqslant #1}							\fi}		
\newcommand{\disjsum}		[1][c]	{
			\ifx p#1 		\amalg 											\else	
			\ifx c#1 		\coprod 										\else	
			\mathrm{Illegal~option} 										\fi\fi}		
			
\newcommand{\fens}			[1]		{\BC{#1}} 										

\newcommand{\Id}			[1]		{\textrm{Id}_{\PM #1}}							
\newcommand{\IdCoef}		[1]		{\fC_{\PM #1}}									
\newcommand{\IntRelB}		[2][]	{\textrm{Int}\BP{#2}}							
\newcommand{\IdMax}			[1]		{\mfM_{\PM #1}}									

\newcommand{\LPlat}			[3][]	{\batonP\BP{#2/#3}} 							
\newcommand{\LNPlat}		[3][]	{\batonQ\BP{#2/#3}_{\PM #1}} 					
\newcommand{\LNPlatS}		[3][]	{\batonQ\BP{#2/#3}_{\PM #1}^{\PM\text{sat}}} 	

\newcommand{\norm}			[2][n]	{
			\ifx n#1 				\left\Vert #2 \right\Vert 				\else
			\ifx o#1 				\left\Vvert #2 \right\Vvert 			\else
			\left\Vert #2 \right\Vert_{\PM #1}								\fi\fi} 

\newcommand{\otimesCSPU}	[1]			{\widehat{\otimes}_{\PM #1}}				
\newcommand{\OuvZa}			[2][]		{\cD_{#1}\BP{#2}} 							

\newcommand{\priv}					{\textbackslash} 								
\newcommand{\PtShilov}		[2][]	{\eta_{\PM #2}}				 					

\newcommand{\reci}			[2][c]	{
			\ifx c#1 		{#2}{}^{\rtimes} 								\else
			\ifx l#1 		{(#2)}^{\rtimes} 								\else
			\mathrm{Illegal~option} 										\fi\fi} 
\newcommand{\rest}			[2]		{#1 | #2} 										

\newcommand{\Sep}			[2][]	{{#2}^{\PM \textbf{sep}}}						
\newcommand{\SpecAnn}		[2][]	{\textrm{Spec}^{#1}\BP{#2}} 					
\newcommand{\SpecB}			[2][]	{\cM\BP{#2}} 									
\newcommand{\SRest}			[3][]	{#2\BC{#3}} 									
\newcommand{\SuppFa}		[1]		{\textrm{Supp}\BP{#1}} 							
\newcommand{\STate}			[3][]	{T_{\PM #2}^{\PM #3}} 							

\newcommand{\Tige}			[3][\Deff]{ 
			\ifx #1\Deff 	#2(#3) 									\else
			\rest{#2}{#1}(#3)		 								\fi}			

\pretitle{\begin{center}\Huge\bfseries}
\posttitle{\par\end{center}\vskip 0.5em}
\preauthor{\begin{center}\Large\ttfamily}
\postauthor{\end{center}}
\predate{\par\large\centering}
\postdate{\par}

\title{Raynaud and Gruson's flattening techniques in Berkovich geometry}
\author{Mattias Ferreira--Filoramo}
\date{}

\newcommand{\Hmm}[1]{\textcolor{red}{${}^\ast$}\leavevmode{\marginpar{\tiny%
            $\hbox to 0mm{\hspace*{-0.5mm}$\leftarrow$\hss}%
            \vcenter{\vrule depth 0.1mm height 0.1mm width \the\marginparwidth}%
            \hbox to 0mm{\hss$\rightarrow$\hspace*{-0.5mm}}$\\\relax\raggedright #1}}} 

\begin{document} \maketitle

\begin{abstract}
\noindent Raynaud and Gruson introduced, in their pioneering paper \cite{RG71}, a flattening procedure for coherent sheaf in the schemes context that rely on sophisticated techniques of algebra and geometry. These methods were later used by Raynaud, Bosch and Lütkebohmert in formal geometry, and more recently Ducros has developed analogous tools for Berkovich spaces. However, it should be noted that the expression "flattening by blowing up" is, in this context, a little abusive because the flattening process actually requires a succession of blow up and quasi-étale morphisms. The aim of this article is to develop a local flattening technique in Berkovich geometry, in the boundaryless case, which only involves blow up.
\end{abstract}

\newpage \tableofcontents

\newpage \section{Presentation of flattening results}

In \cite{RG71}, Raynaud and Gruson introduced a flattening procedure for families of varieties based on sophisticated techniques from algebraic geometry. More specifically, this procedure was introduced for coherent sheaves on a scheme and was subsequently used by Raynaud, Bosch and Lütkebohmert to develop rigid analytical geometry using formal models. In this article, we continue to develop, in the context of Berkovich spaces, the flattening techniques inspired by those of Raynaud and Gruson and developed by Antoine Ducros in \cite{Duc20} and \cite{Duc23}.

\subsection{The flattening results of Raynaud and Gruson from \cite{RG71}}

Before presenting the flattening results in non-Archimedean geometry, let us first outline a simplified version of those of Raynaud and Gruson. Let $X$ and $Y$ be two Noetherian schemes, $f$ a finite-type morphism from $Y$ to $X$, $\cF$ a coherent $\cOf{Y}$-module and $U$ a quasi-compact open set in $X$ above which $\cF$ is $X$-flat, Theorem (5.2.2) of \cite{RG71} guarantees the existence of a closed subscheme $T$ of $X$, with support in $X \priv U$, such that the strict transform of is flat over the blow-up of $X$ along $T$.\par
Let $\tilde{X}$ be the blow-up of $X$ along $T$, it should be noted that the image of the blow-up morphism $\pi$ from $\tilde{X}$ to $X$ is equal to the adherence of $U$ in $X$. It is therefore when $U$ is dense that this theorem of Raynaud and Gruson is most powerful. In the opposite case, where $U$ is empty, the scheme $\tilde{X}$ is also empty and the theorem is meaningless. Thus, these flattening techniques do not create flatness out of nothing: they merely propagate the flatness that is already present to some extent.\par\medskip

Raynaud and Gruson’s result also allows us to control the dimension of the relative non-flatness locus of $\cF$. This leads to the concept of a flat module in controled dimensions:

\begin{deff}[Flat module in controled dimensions]
Let $X$ and $Y$ be two schemes, $f$ a finitely presented morphism from $Y$ to $X$, $\cF$ a $\cOf{Y}$-module of finite type and $n$ an integer, the $\cOf{Y}$-module $\cF$ is said to be $X$-flat in dimensions greater than $n$ if there exists a retro-compact open set $V$ in $Y$ such that $\rest{\cF}{V}$ is a $\cOf{V}$-module of finite presentation that is $X$-flat and suche that the relative dimension $\dimKr{(Y \priv V)/X}$ of $Y \priv V$ over $X$ is strictly less than $n$.
\end{deff}

If we know of an upper bound on the dimension of the locus of relative non-flatness of $\cF$, then Raynaud and Gruson’s result states that the strict transform $\tilde{\cF}$ of $\cF$ with respect to $\pi$ satisfies the same upper bound on the locus of relative non-flatness.

\begin{theo}[Theorem (5.2.2) of \cite{RG71}]
Let $X$ be a quasi-compact and quasi-separated scheme, $U$ a quasi-compact open set in $X$, $(Y,f)$ an $X$-finite presentation scheme, $\cF$ a $\cOf{Y}$-module of finite type and $n$ an integer, if $\rest{\cF{}}{\reci{f}\fens{U}}$ is $U$-flat in dimensions greater than $n$, then there exists a $U$-admissible blow-up $(\tilde{X},\pi)$ such that the strict transform of $\cF$ by $\pi$ is $\tilde{X}$-flat in dimensions greater than $n$.
\end{theo}

\[ \xymatrix{
& \rest{\cF{}}{\reci{f}\fens{U}} \ar@{..}[dl] && \cF \ar@{..}[dl] && \tilde{\cF} \ar@{..}[dl] \\
\reci{f}\fens{U} \ar[rr] \ar[d]_{\rest{f}{\reci{f}\fens{U}}} && Y \ar[d]_f && Y \times_X \tilde{X} \ar[d] \ar[ll] \\
U \ar[rr] && X && \ar[ll]^\pi \tilde{X}
} \]

Let $\tilde{Y}$ denote the scheme $Y \times_X \tilde{X}$, $\batonP$ the open set in $Y$ consisting of the points at which $\cF$ is $X$-flat, and $\batonQ$ the complement of $\batonP$ in $Y$, noting that $T$ is set-theoretically equal to $\adh{f\fens{\batonQ}}$, we can reformulate some properties of $\pi$ in a way that may seem more obscure, but which will be useful for stating our main result.
\begin{enumerate}[itemsep=1mm]
\item The quotient of the inverse image of $\cF$ on $\tilde{Y}$ by its subshaef of sections set-theoretically supported in $\tilde{Y} \times_X T$ is $\tilde{X}-$flat.
\item The open set $\tilde{Y} \priv (\tilde{Y} \times_X T)$ of $\tilde{Y}$ lies above $\batonP$.
\item The closed set $\tilde{Y} \times_X T$ of $\tilde{Y}$ lies set-theoretically above $\reci{f}\fens{\adh{f\fens{\batonQ}}}$.
\item The image of $\tilde{Y}$ in $Y$ contains $Y \priv \reci{f}\fens{\adh{f\fens{\batonQ}}}$.
\item For any irreducible component $G$ of $\adh{f\fens{\batonQ}}$ that is not an irreducible component of $X$, the image of $\tilde{Y} \times_X T$ in $Y$ contains $\reci{f}\fens{G}$ set-wise.
\end{enumerate}

\subsection{The flattening results in non-Archimedean geometry from \cite{Duc20} and \cite{Duc23}}

Starting from a flat coherent sheaf in the sense of rigid geometry, Raynaud, Bosch and Lütkebohmert constructed a flat formal model using decompositions centred on the special fibre. Recently, Ducros has developed similar tools in non-Archimedean geometry, not within the rigid framework, but within that of Berkovich spaces. In \cite{Duc18}, Ducros conducts a systematic study of flatness in this context and introduces a notion of "dévissage" in the style of Raynaud and Gruson. In \cite{Duc20}, he develops the concept of flattening via blow-ups and provides some applications to the study of images of morphisms between compact analytic spaces.\par
Although his methods are similar to those originally developed by Raynaud and Gruson, the result is much more difficult to achieve in the context of Berkovich spaces. Indeed, the term "flattening by blow-ups" can, in the context of Berkovich spaces, be somewhat misleading: although these techniques are directly inspired by those of Raynaud and Gruson, they give rise to statements that are considerably more complicated than their schematic originals. In \cite{Duc20}, Theorem (6.6) allows a coherent module to be flattened solely by strict transformation after a change of basis along a sequence of blow-ups and quasi-étale morphisms, which it is probably not possible to avoid according to (7.12) of \cite{Duc20}.\par
Although the proof of Theorem (6.6) in \cite{Duc20} broadly follows the strategy of Raynaud and Gruson, which is based on the concept of dévissage, its application in the context of Berkovich spaces encounters a number of significant difficulties, which have been resolved in some of the author’s earlier works and in \cite{Duc20}.

\subsection{Main results}

To begin with, we continue to develop, within the context of Berkovich spaces, flattening techniques based on blow-ups, inspired by those introduced by Raynaud and Gruson in their seminal work \cite{RG71}, and further developed by Antoine Ducros in \cite{Duc20} and \cite{Duc23}. We thus develop a local flattening technique involving only blow-ups and immersions of affinoid domains — these are local results — for spaces without a boundary. More specifically, under the assumption of the absence of relative boundaries, we obtain a local version of the flattening result without the quasi-étale morphisms present in \cite{Duc20}. This constitutes the central theorem of this work.\par
The concept of a boundary of a non-Archimedean analytic space is intrinsically an analytic concept, with no equivalent in algebraic geometry. For example, a closed disc has a boundary, whereas an open disc does not. Since every algebraic variety lacks a boundary, and given that Raynaud and Gruson’s original results on flattening were established in the schematic framework, it was reasonable to expect analogous conclusions (in particular, the elimination of quasi-étale morphisms) for analytic spaces without a boundary.\par
Our new result on local flattening provides a much more flexible and manageable variant of the general statement, both in terms of the flattening itself and its various possible applications, under assumptions that remain very reasonable. In practice, as the boundary is a relative concept, it is very common to work with boundaryless morphisms, or at least overconvergent morphisms—that is, those which admit a boundaryless extension. This is notably the case for quasi-smooth morphisms.\par\bigskip

Next, we apply our flattening theorem to morphisms between analytic spaces, in order to derive some applications in line with those given in \cite{Duc20} and \cite{Duc23}. We use it to locally improve the structure of a pointed morphism $f$ between reduced pointed analytic spaces $(Y,y)$ and $(X,x)$, such that the source $Y$ is equidimensional of dimension $n$, by requiring the existence of an integer $d$ such that $f$ is generically of dimension $d$. These assumptions hold in a canonical manner whenever the source is irreducible, which is the reference configuration for this study.\par
We show that there exists a morphism $\varphi$ from a $K$-analytic space $Z$ to $X$ which is a composition of blow-ups and inclusions of analytic domains, a divisor $S$ of $Z$ and an analytic neighbourhood $U$ of $y$ in $Y$, such that, denoting by $V$ the $K -$analytic space $U \times_X Z$, we have:
\begin{enumerate}
\item The strict transformation of $V$ by $\varphi$ is surjective onto $U$.
\item The morphism of the strict transformation of $V$ by $\varphi$ into $Z$ factors as a flat and surjective morphism onto a reduced closed analytic subspace $F$ of $Z$ of pure dimension $n-d$.
\item The image of the exceptional divisor of $\varphi$ in $X$ is contained in that of $U$ and has a dimension strictly less than $n-d$.
\end{enumerate}\par\bigskip

We then turn our attention to the local structure of the image of a pointed morphism $f$ between compact pointed analytic spaces $(Y,y)$ and $(X,x)$. Until the recent developments in Section 5 of \cite{Duc23}, the geometry of this image had, to our knowledge, only been elucidated under assumptions of flatness or properness. For a flat morphism, the image constitutes a compact analytic domain of $X$, in accordance with (5.11) of \cite{BL93} for the strict setting and with (9.2.1) of \cite{Duc18} for the general case. In the event that $f$ is proper, the image is a Zariski closed set of $X$; on this point, we refer the reader to Section 1.3.23 of \cite{Duc18}.\par
To describe the image of $f$ locally, we introduce the concept of an $m$-cell on $X$, defined by recurrence on the integer $m$. A $0$-cell of $X$ consists of a finite, non-empty set of rigid points of $X$. For $m$ a non-zero integer, an $m$-cell of $X$ is a non-empty, locally closed subset $C$ of $X$ satisfying the following properties:
\begin{enumerate}
\item There exists a finite family of closed analytic domains $(V_i)_{i \in I}$ of $X \priv (\adh{C} \priv C)$ and a non-empty, purely $n$-dimensional, Zariski closed set $Z$ of the disjoint sum of  $(V_i)_{i \in I}$, whose image in $X$ is $C$.
\item There exists a finite family of integers $(n_j)_{j \in J}$ strictly less than $n$ and a family $(C_j)_{j \in J}$ of subsets of $X$ such that, for every $j$ in $J$, the set $C_j$ is an $n_j$-cell and $(\adh{C} \priv C)$ is the union of $(\adh{C_j})_{j \in J}$.
\end{enumerate}
We show that there exists an analytic neighbourhood $U$ of $y$ in $Y$, an analytic neighbourhood $V$ of $x$ in $X$, and an $(n-d)$-cell $C$ of $V$, such that:
\begin{enumerate}
\item The image of $U$ under $f$ is equal to $\adh{C}$.
\item Let $\Abhy{U}$ denote the set of Abhyankar points of $U$, the set $f\fens{\Abhy{U}}$ is contained in $C$.
\end{enumerate}
We view this decomposition as a kind of local, non-Archimedean Chevalley theorem.\par\bigskip

It should be noted, however, that reasonable stratifications of images of analytic morphisms were constructed by Lipshitz in a completely different language and using entirely different methods, his framework being the theory of models of algebraically closed valued fields in a language enriched by a class of "reasonable" analytic functions. We refer in particular to (5.2) in \cite{Lip93} and are unaware of any precise links between his structural theorems and our own.\par\bigskip

\textbf{Use of Large Language Models.} During the preparation of this article, the author used Deepl only as interactive tools to assist with translation from french to english. The author assume full responsibility for the contents of the article.

\newpage \section{Reminders, conventions and notations}

The results presented in this section are mainly drawn from the seminal work of Vladimir Berkovich and various works by Antoine Ducros. We review here some prerequisites of the theory of analytic spaces in the sense of Berkovich, as developed in \cite{Ber90} and \cite{Ber93}, the reading of which may be usefully guided by the notes \cite{Jon19}. For the study of local theory — non-Archimedean analysis, seminorms and commutative Banach algebras — the standard reference remains \cite{BGR}.

\subsection*{General conventions}

Let $X$ and $Y$ be two sets, $f$ a function from $X$ to $Y$, $U$ a subset of $X$ and $V$ a subset of $Y$. We denote $f\fens{U}$ the direct image of $U$ under $f$ and $\reci{f}\fens{V}$ the reciprocal image of $V$ under $f$. Furthermore, for any element $v$ of $V$, we denote $\reci{f}(v)$ as the inverse image of $v$ under $f$.\par\medskip

Let $a$ and $b$ be two integers, we denote $\lseg a,b \rseg$ the set of integers between $a$ and $b$.

\subsection*{Algebraic conventions}

\noindent In this article, unless otherwise specified, all groups are assumed to be commutative, and their composition laws are denoted additively. Furthermore, all rings are assumed to be commutative and unitary. Finally, all algebras are assumed to be unitary and associative, all algebra homomorphisms are assumed to be unitary and by a subalgebra of an algebra $A$, we mean a subalgebra containing the identity element of $A$.\par\medskip

Let $A$ be a ring, we denote by $A^+$ the additive group underlying $A$, by $A^\times$ the multiplicative group of invertible elements of $A$ and by $\SpecAnn{A}$ the spectrum of $A$.\par\medskip

Let $A$ be a Noetherian local ring and $M$ a finite-type $A$-module, we denote $\Coprofl[A]{M}$ by the codepth of $M$.

\subsection*{Topological conventions}

\noindent In this article, we follow Bourbaki’s conventions regarding topological spaces \cite{Bou07}.\par\medskip

Let $X$ be a topological space and $E$ a subset of $X$, we denote $\adh[X]{E}$ as the closure of $E$ in $X$.\par\medskip

A topological space is said to be quasi-compact if every open cover of $X$ has a finite subcover. It is said to be compact if it is both quasi-compact and separated.\par\medskip

Let $X$ and $Y$ be two topological spaces and $f$ a continuous map from $X$ to $Y$. We say that $f$ is topologically proper if $f$ is universally closed, that is, if, for every topological space $Z$, the map $f \times \Id{Z}$ is closed. The reader should note that this definition does not require $f$ to be separated, although this distinguishes it from the concept of a proper morphism in algebraic geometry.\par\medskip

Let $X$ and $Y$ be two topological spaces and $f$ a continuous map from $X$ to $Y$. For $f$ to be proper, it is necessary and sufficient that $f$ be closed and that, for every $y$ in $Y$, the set $\reci{f}(y)$ be quasi-compact. It follows that, for any proper map, the inverse image of any quasi-compact set is quasi-compact by (I, § 10, no. 2, Prop. 6) of \cite{Bou07}.\par
Assusme that $X$ is separated and that $Y$ is locally compact. For $f$ to be proper, it is necessary and sufficient that the inverse image of every compact set under $f$ is compact, by (I, \S 10, \no 3, Prop. 7) of \cite{Bou07}. Furthermore, if $f$ is proper, then $X$ is locally compact.\par\medskip

Let $X$ be a topological space and $x$ a point in $X$. We denote $\dimKr{X}$ as the Krull dimension of $X$ and $\dimKr[x]{X}$ as the Krull dimension of $X$ at the point $x$.\par\medskip

Let $(X,\cA)$ be a locally ringed space and $x$ a point in $X$. We denote by $\cA_x$ the fiber of $\cA$ at $x$, by $\IdMax{x}$ the maximal ideal of $\cA_x$, and by $\CrpRes{x}$ the residue field $\cA_x/\IdMax{x}$ of $\cA_x$. Furthermore, for any global section $f$ of $\cA$, we denote by $f_x$ the image of $f$ in $\cA_x$.

\subsection*{Uniforms conventions}

\noindent In this article, we follow Bourbaki’s conventions regarding uniform spaces \cite{Bou07}.\par\medskip

When a topological group is considered as a uniform space, it always refers—unless otherwise specified—to the uniform structure associated with its topology. Similarly, when we consider the uniform structure of a topological ring, we are always referring to the uniform structure of its additive group. In particular, a ring is said to be complete if its additive group is complete.\par\medskip

Let $G$ be a topological group, we denote $\Sep{G}$ as the separated group of $G$ and $\CoSep{G}$ as the completed separated group of $G$.

\subsection*{Analytic spaces conventions}

For the remainder of this section, we fix a complete ultrametric field $K$, whose valuation may be trivial. We will work with the concept of a $K$-analytic space in the sense of Berkovich and assume familiarity with the fundamentals of the theory, as set forth, for example, in the seminal texts \cite{Ber90} and \cite{Ber93}. Our definition of a $K$-analytic space will be that given in \cite{Ber93}, and with this choice, the spaces defined in \cite{Ber90} are the proper $K$-analytic spaces, that is, those in which every point has an affinoid neighborhood.\par\medskip

We also fix a subgroup $\Gamma$ of $\eR^\times_+$ such that $\abs{K^\times}\!{}^{\eQ} \cdot \Gamma$ is not trivial. In other words, the group $\Gamma$ is nontrivial if the absolute value of $K$ is trivial. We will use the concept of a $\Gamma$-strict $K$-analytic space, for which we refer to Section (3.1) of \cite{Duc18}. Intuitively, a $K$-analytic space is $\Gamma$-strict if it can be defined using only real parameters belonging to $\Gamma$. Thus, if $\Gamma$ is trivial, the $\Gamma-$strict $K-$analytic spaces are the strictly $K-$analytic spaces, and if $\Gamma$ is equal to $\eR^\times_+$, then every $K-$analytic space is $\Gamma-$strict.

\subsection*{Basic notations for analytic spaces}

Let $X$ be a $K$-analytic space.\par\medskip

Let $x$ be a point in $X$. We denote by $\CrpResC{x}$ the residue field completed at $x$ and by $\CrpResCR{x}$ the reduction of $\CrpResC{x}$.\par\medskip

Let $L$ be a complete extension of $K$, we denote by $X_L$ the $L$-analytic space obtained by extending the scalars from $K$ to $L$. There exists a $K$-morphism from $X_L$ to $X$ that is surjective by (0.5) of \cite{Duc07}. When $X$ is a $K$-affinoid space of algebra $A$, this morphism is induced by the morphism from $A$ to $A \otimesCSPU{K} L$, which is faithfully flat by (2.1.2) of \cite{Ber93}.\par\medskip

Let $Y$ be a $K$-analytic space and $f$ a morphism from $Y$ to $X$, we denote by $\reci{f}(x)$ or $Y_x$ the fiber of $f$ at $x$. This is a $\CrpResC{x}$-analytic space.\par\medskip

\noindent An analytic space without a base field is a pair $(X,L)$ where $L$ is a complete ultrametric field and $X$ is an $L$-analytic space. Let $(X,L)$ and $(Y,E)$ be two analytic spaces, a morphism from $(Y,E)$ to $(X,L)$ consists of an isometric embedding of $L$ into $E$ and a morphism of $E$-analytic spaces from $Y$ to $X_E$.

\begin{att}
Let $(X,L)$ be an analytic space; the $L$-analytic space $X$ is not identical to $(X,L)$, but for the sake of readability and clarity, we commonly identify them. We will then say “let $X$ be an analytic space” and refer to $L$ as the domain of definition of $X$. Similarly, for morphisms between analytic spaces, we will not mention the domains of definition or the associated isometric embedding.
\end{att}

\textbf{GAGA} Let $A$ be a $K$-affinoid algebra and $X$ be the $K$-affinoid space $\SpecB{A}$. We will denote $\Algbr{X}$ the affine scheme given by the ring $A$, it is equipped with a morphism of locally ringed spaces from $X$ to $\Algbr{X}$.

\newpage \section{Local factorization of morphisms}

\begin{att}
Throughout this section, $K$ denotes a non-trivially valued field.
\end{att}

\begin{deff}[Noether normalization]
Let $X$ be a $K$-analytic space, a Noether normalization of $X$ is a finite, surjective morphism from $X$ to $\eD^d_K(0,1)$.
\end{deff}

\begin{nott}
We denote $\STate{n}{K}$ the normed $K$-algebra\footnote{We then set $\STate{0}{K} = K$.} $\SRest{K}{X_1, \ldots, X_n}$.
\end{nott}

\begin{lem}[Perturbation of a Noetherian normalization]\label{LEM - deformation noether affinoide}
Let $A$ be a strictly $K$-affinoid algebra and $\pi$ a contractive, admissible, injective, finite homomorphism from $\STate{d}{K}$ to $A$. There exists a nonzero positive real number $\varepsilon$ such that, for any family $(g_i)_{1 \leqslant i \leqslant d}$ of elements of $A$ satisfying that, for every $i$ between $1$ and $d$, we have:
\[ \norm[A]{ \pi(X_i) - g_i } \leqslant \varepsilon \]
the ring homomorphism $\psi$ from $\STate{d}{K}$ to $A$ defined by the family $(g_i)_{1 \leqslant i \leqslant d}$, and coinciding with $\pi$ in $A$, exists and is a finite, injective, contractive, lipschitz homomorphism.
\end{lem}
\begin{dem}
The existence of a nonzero positive real number $\varepsilon$ such that, for any family $(g_i)_{1 \leqslant i \leqslant d}$ of elements of $A$ satisfying, for every $i$ between $1$ and $d$, the following:
\[ \norm[A]{ \pi(X_i) - g_i } \leqslant \varepsilon \]
The ring homomorphism $\psi$ from $\STate{d}{K}$ to $A$, defined by the family $(g_i)_{1 \leqslant i \leqslant d}$ and coinciding with $\pi$ in $A$, exists and is a contractive homomorphism, as follows from (7.3.4/3) of \cite{BGR}. It remains to show that $\psi$ is a finite ring homomorphism from $\STate{d}{K}$ to $A$, and to do so, note that the residual operator $\AnnRes[i]{\psi}$ is equal to $\AnnRes[i]{\pi}$. It follows from property (2.5.2) of \cite{Ber90} that the morphism between $K$-affinoid spaces induced by $\psi$ is boundaryless, and is therefore finite by (2.5.13) of \cite{Ber90}
\end{dem}

\begin{lem}[Induced Noetherian normalization]\label{LEM - deformation noether}
Let $X$ be a strictly $K$-affinoid space, $Y$ a strictly $K$-affinoid $X$-space and $x$ a point in $X$. There exists an integer $d$, an affinoid neighborhood $U$ of $x$ in $X$ and a quasi-finite $U$-morphism $\Psi$ from $Y \times_X U$ to $\eA^d_U$ that induces a Noetherian normalization of the $\CrpResC{x}$-analytic space $Y_x$.
\end{lem}
\begin{dem}
Since the fiber $Y_x$ is strictly affinoid, there exists, by Noether’s normalization theorem 6.1.2/2 in \cite{BGR}, an integer $d$ and a finite surjective morphism $\psi$ from $Y_x$ to the unit polydisk of dimension $d$ on $\CrpResC{x}$. It follows, even if we have to restrict $X$ to the neighborhood of $x$ and deform $\psi$, that there exists a quasi-finite $X$-morphism $\Psi$ from $Y$ to $\eA^d_X$ that induces $\psi$ on $Y_x$ by lemma \ref{LEM - deformation noether affinoide}.
\end{dem}

\begin{prop}[Characterization of the boundary by reduction]\label{Carc-bord-red}
Let $X$ and $Y$ be two $K$-affinoid spaces, $f$ a morphism from $Y$ to $X$, and $y$ a point in $Y$, the point $y$ belong to the interiror of $f$, if and only if the closure of $\tilde{y}$ in $\widetilde{Y}$ is a finite $\widetilde{X}$-scheme.
\end{prop}
\begin{dem}
Let $\AnnRes[i]{A}$ be the ring of global sections of $\widetilde{X}$ and $\AnnRes[i]{B}$ be the ring of global sections of $\widetilde{Y}$, for $y$ to belong to $\IntRelB{f}$, it is necessary and sufficient that the image of $\AnnRes[i]{B}$ in $\CrpRes{\tilde{y}}$ be integral over $\AnnRes[i]{A}$ by the fourth point of property (2.5.2) in \cite{Ber90}. Now, since the image of $\AnnRes[i]{B}$ in $\CrpRes{\tilde{y}}$ defines a closed subscheme of $\widetilde{Y}$ corresponding to the adherence of $\tilde{y}$ in $\widetilde{Y}$, we obtain the desired equivalence.
\end{dem}

\begin{prop}\label{PROP - factorisation fini et lisse connexe}
Let $(X,x)$ and $(Y,y)$ be two good, strict $K$-analytic pointed spaces, $f$ a pointed morphism from $(Y,y)$ to $(X,x)$ such that $y$ lies in the interior $f$ and $d$ the dimension of $f$ at $y$, there exists a strict rational affinoid neighborhood $V$ of $y$ in $Y$, a strict affinoid neighborhood $W$ of $x$ in $X$ and a $d$-polyradius $r$ with coefficients in $\abs{K^\times}\! {}^{\eQ}$ such that $\rest{f}{V}$ factors through a finite morphism from $V$ to $\eD^d_W(0,r)$.
\end{prop}
\begin{dem}
Since $K$ is not trivially valued and the result is local on $X$ and $Y$, we may assume that $X$ and $Y$ are strictly $K$-affinoid. Let $\Delta$ denote the unit polydisk of dimension $d$ on $\CrpResC{x}$, even if this requires restricting $X$ to the neighborhood of $x$, there exists a quasi-finite $X$-morphism $\Psi$ from $Y$ to $\eA^d_X$ whose induced morphism $\psi$ on $Y_x$ is a Noetherian normalization of $Y_x$ by the lemma \ref{LEM - deformation noether} :
\[ \xymatrix{
Y \ar[d]^\Psi && \ar[ll] Y_x \cartesien \ar[d]^\psi \\
\eA^d_X && \ar[ll] \Delta } \]
The reduction of $\Delta$ is the $d$-dimensional affine space over $\CrpResCR{x}$ and the morphism $\tilde{\psi}$ from $\widetilde{Y_x}$ to $\widetilde{\Delta}$, the reduction of $\psi$, is finite by 6.3.5/1 of \cite{BGR}.
\[ \xymatrix{
Y \ar[d]^\Psi && \ar[ll] Y_x \cartesien \ar[d]^\psi \ar@{.>}@/^1pc/[rrrr]^{\text{red}} && && \widetilde{Y_x} \ar[d]^{\tilde{\psi}} \\
\eA^d_X && \ar[ll] \Delta && && \widetilde{\Delta} } \]
Now, since the point $y$ lies inside $f$, the reduction $\tilde{y}$ is a closed point of $\widetilde{Y_x}$ by property \ref{Carc-bord-red}. Furthermore, since $\tilde{\psi}$ is closed, the point $\tilde{\psi}(\tilde{y})$ is closed. Let $\tilde{\pi}$ be the endomorphism of the $d$-dimensional affine space on $\CrpResCR{x}$ given by the minimal polynomials of the coordinates of $\tilde{\psi}(\tilde{y})$. The morphism $\tilde{\pi}$ is finite, and we have:
\[ (\tilde{\pi} \circ \tilde{\psi})(\tilde{y}) = 0 \]
Even if we restrict $X$ to the neighborhood of $x$, there exists an endomorphism $\Pi$ of $\eA^d_X$ that induces an endomorphism $\pi$ of the unit polydisk of dimension $d$ onto $\CrpResC{x}$, whose reduction yields the morphism $\tilde{\pi}$. Indeed, since $\kappa(x)$ is dense in $\CrpResC{x}$, it is possible to choose $\pi$ such that its reduction is $\tilde{\pi}$ and it is defined on $\kappa(x)$.
\[ \xymatrix{
Y \ar[d]^\Psi && \ar[ll] Y_x \cartesien \ar[d]^\psi \ar@{.>}@/^1pc/[rrrr]^{\text{red}} && && \widetilde{Y_x} \ar[d]^{\tilde{\psi}} \\
\eA^d_X \ar[d]^\Pi && \ar[ll] \Delta \cartesien \ar[d]^{\pi} && && \widetilde{\Delta} \ar[d]^{\tilde{\pi}} \\
\eA^d_X && \ar[ll] \Delta && && \widetilde{\Delta} } \]
If we replace $\psi$ with $\pi \circ \psi$ and $\Psi$ with $\Pi \circ \Psi$, we can assume that $\tilde{\psi}(\tilde{y})$ is zero. Let $\Delta\!{}^\circ$ denote the open unit polydisk of dimension $d$ on $\CrpResC{x}$. Since the morphism $\tilde{\psi}$ is finite, the fiber of $\tilde{y}$ is finite, and since $\tilde{y}$ is closed in $\widetilde{Y_x}$, the set:
\[ C = \reci{\text{red}}(\tilde{y}) \]
contains $y$ and is a closed open set in $Y_x \times_\Delta \Delta\!{}^\circ$:
\[ \xymatrix{
Y \ar[d]^\Psi && \ar[ll] Y_x \ar[d]^\psi && C \ar[ll] \ar[d]^\psi && \\
\eA^d_X && \ar[ll] \Delta && \Delta\!{}^\circ \ar[ll] } \]
Since every point in $C$ maps to $\tilde{y}$ and since the point $y$ lies in the interior of $f$, the point $y$ lies in the interior of $\psi$ by (1.5.5) of \cite{Ber93}, and the set $C$ avoids the boundary of $f$ by the property \ref{Carc-bord-red}. Let $t$ denote the image of $y$ under $\psi$; since the reduction of $t$ is zero, the point $t$ belongs to the open unit polydisk of dimension $d$ on $\CrpResC{x}$. There then exists a $d$-polyradius $r$ with coefficients in $\abs{K^\times}\!{}^{\eQ}$, whose components are strictly less than $1$, such that $t$ belongs to $\eD^d_{\CrpResC{x}}(0,r)$. Let $U$ be the set:
\[ U = C \times_\Delta \eD^d_{\CrpResC{x}}(0,r) \]
The set $U$ is a closed open set in $Y_x \times_\Delta \eD^d_{\CrpResC{x}}(0,r)$ and contains $y$:
\[ \xymatrix{
Y \ar[d]^\Psi && \ar[ll] Y_x \ar[d]^\psi && C \ar[ll] \ar[d]^\psi && U \cartesien \ar[ll] \ar[d]^\psi \\
\eA^d_X && \ar[ll] \Delta && \Delta\!{}^\circ \ar[ll] && \eD^d_{\CrpResC{x}}(0,r) \ar[ll] } \]
Since $U$ is a closed open set, there exists an idempotent $e$ in $\cO(Y_x \times_\Delta \eD^d_{\CrpResC{x}}(0,r))$ such that:
\[ U = \{ \abs{e} = 0 \} \]
However, since $Y_x \times_\Delta \eD^d_{\CrpResC{x}}(0,r)$ is a Weierstrass domain of $Y_x$ by 2.2.2 of \cite{Ber90}, the image:
\[ \left| \begin{array}{c c c}
\cO(Y_x) & \to & \cO(Y_x \times_\Delta \eD^d_{\CrpResC{x}}(0,r))
\end{array} \right. \]
is dense by 2.2.10 of \cite{Ber90}. Let $\eta$ be a real number strictly positive and strictly less than $1/2$, and let $\varepsilon$ be a real number strictly positive and less than $\sqrt{\eta}$. By restricting $X$ to the neighborhood of $x$ and by the density of $\kappa(x)$ in $\CrpResC{x}$, there exists a function $\gamma$ in $\cO(Y)$ such that its restriction $\gamma_0$ to $Y_x$ satisfies:
\[ \norm{\gamma_0 - e} \leqslant \varepsilon, \quad \norm{\gamma_0} = \norm{e} \leqslant 1, \quad \norm{1 - \gamma_0} = \norm{1 - e} \leqslant 1 \quad\et\quad \norm{\gamma_0^2 - \gamma_0} < \eta \]
Since $\{ \abs{\gamma^2 - \gamma} \geqslant \eta \}$ is compact, even if we have to restrict $X$ to the neighborhood of $x$, we can assume that the following holds:
\[ Y = \{ \abs{\gamma^2 - \gamma} < \eta \} \]
Let $z$ be a point in $Y$. We have $\abs{\gamma^2(z) - \gamma(z)} < \eta$, and given the triangle inequality and the fact that $\eta$ is strictly less than $1$, we conclude that the following two conditions :
\[ \abs{\gamma(z)} \leqslant \sqrt{\eta} \quad\et\quad \abs{1 - \gamma(z)} \leqslant \sqrt{\eta} \]
are mutually exclusive. Setting:
\[ V_0 = \{ \abs{\gamma} \leqslant \sqrt{\eta} \} \quad\et\quad V_1 = \{ \abs{1 - \gamma} \leqslant \sqrt{\eta} \} \]
We have the following breakdown:
\[ Y = V_0 \sqcup V_1 \]
The set $V_0$ is therefore a closed open set of $Y$ containing $U$. Let $V$ denote the set:
\[ V = V_0 \times_{\eA^d_X} \eD^d_X(0,r) \]
The set $V$ is a closed open set in $Y \times_{\eA^d_X} \eD^d_X(0,r)$. Since $\varepsilon$ is less than $\sqrt{\eta}$ and $\norm{\gamma_0 - e}$ is less than $\varepsilon$, the following relations :
\[ \abs{e} \leqslant \sqrt{\eta} \quad\et\quad \abs{\gamma_0} \leqslant \sqrt{\eta} \]
as well as the relations:
\[ \abs{1 - e} \leqslant \sqrt{\eta} \quad\et\quad  \abs{1 - \gamma_0} \leqslant \sqrt{\eta} \]
are equivalent. We can then deduce the following equalities:
\[ V \times_{\eD^d_X(0,r)} \eD^d_{\CrpResC{x}}(0,r) = \{ \abs{e} \leqslant \sqrt{\eta} \} \cap ( Y_x \times_\Delta \eD^d_{\CrpResC{x}}(0,r) ) \]
and 
\[ V_1 \times_{\eD^d_X(0,r)} \eD^d_{\CrpResC{x}}(0,r) = \{\abs{1-e} \leqslant \sqrt{\eta}\} \cap (Y_x \times_\Delta \eD^d_{\CrpResC{x}}(0,r)) \]
Furthermore, since we have the disjoint union:
\[ Y_x \times_\Delta \eD^d_{\CrpResC{x}}(0,r) = ( V \times_{\eD^d_X(0,r)} \eD^d_{\CrpResC{x}}(0,r)) \sqcup ( V_1 \times_{\eD^d_X(0,r)} \eD^d_{\CrpResC{x}}(0,r) ) \]
and since $V_1 \times_{\eD^d_X(0,r)} \eD^d_{\CrpResC{x}}(0,r)$ contains the complement of $U$ in $Y_x \times_\Delta \eD^d_{\CrpResC{x}}(0,r)$, we have:
\[ U = V \times_{\eD^d_X(0,r)} \eD^d_{\CrpResC{x}}(0,r) \]
and the commutative diagram:
\[ \xymatrix{
Y \ar[d]^\Psi && \ar[ll] V \ar[d]^\Psi && U \cartesien \ar[ll] \ar[d]^\psi \\
\eA^d_X && \ar[ll] \eD^d_X(0,r) && \ar[ll] \eD^d_{\CrpResC{x}}(0,r) } \]
Let $B$ denote the set $\BorRelB{f} \cap V$. Since $V$ is compact, the set $B$ is compact and does not intersect $U$. Furthermore, the sets $f(V)$ and $f(B)$ are compact, and $x$ does not belong to $f(B)$. By restricting $X$ around $x$, we can therefore assume that $V$ does not intersect the boundary of $f$, and the morphism from $V$ to $\eD^d_X(0,r)$ is boundaryless. The morphism that $\Psi$ induces on $V$ is quasi-finite and boundaryless; therefore, by 2.5.13 of \cite{Ber90}, it is finite, and we have the desired result.
\end{dem}

\begin{lem}[Pullback of rational domains]\label{LEM - Relevement de domaines rationnels}
Let $(X,x)$ be a strict $K$-analytic, pointed good space and $Z$ be a Zariski closed set of $X$ containing\footnote{In particular, the space $(Z,x)$ is good and strictly $K$-analytic.} $x$ and $V$ a strict rational affinoid neighborhood of $x$ in $Z$, there exists a strict rational affinoid neighborhood $W$ of $x$ in $X$ such that $V = Z \cap W$.
\end{lem}
\begin{dem}
Let $(f_i)_{1 \leqslant i \leqslant n}$ and $g$ be analytic functions on $Z$ with no common zeros such that:
\[ V = \{ \abs{f_1} \leqslant \abs{g}, \quad\ldots,\quad \abs{f_n} \leqslant \abs{g} \} \]
Since the functions $(f_i)_{1 \leqslant i \leqslant n}$ and $g$ have no common zero, the function $g$ does not vanish on $V$ and is therefore bounded below by a real number $c$ belonging to $\abs{K^\times}\!{}^{\eQ}$. We can then add the following inequality:
\[ \abs{1} \leqslant c^{-1} \abs{g} \]
to the system describing $V$. For every $i$ between $1$ and $n$, let $f_i$ be an analytic function $\varphi_i$ on $X$, and let $g$ be an analytic function $\psi$ on $X$. The rational domain $W$ of $X$ is defined by:
\[ W = \{ \abs{\varphi_1} \leqslant \abs{\psi}, \quad\ldots,\quad \abs{\varphi_n} \leqslant \abs{\psi} \quad\et\quad \abs{1} \leqslant c^{-1} \abs{\psi} \} \]
satisfies the desired property.
\end{dem}

\begin{prop}\label{PROP - factorisation fini et lisse connexe sur un ferme}
Let $(Y,y)$ and $(X,x)$ be two strict $K$-analytic pointed spaces, $Z$ be a Zariski closed set in $Y$ containing\footnote{In particular, the space $(Z,y)$ is good and strict.} $y$, $f$ be a pointed morphism from $(Z,y)$ to $(X,x)$ such that $y$ lies in the interior of $f$ and $d$ be the dimension of $f$ at $y$. There exists a strict rational affinoid neighborhood $V$ of $y$ in $Y$, a strict affinoid neighborhood $W$ of $x$ in $X$, and a $d$-polyradius $r$ with coefficients in $\abs{K^\times}\! {}^{\eQ}$ such that $\rest{f}{(V \times_Y Z)}$ factors through a finite morphism from $V \times_Y Z$ to $\eD^d_W(0,r)$.
\end{prop}
\begin{dem}
The result follows from property \ref{PROP - factorisation fini et lisse connexe} and lemma \ref{LEM - Relevement de domaines rationnels}.
\end{dem}

\newpage \section{Control of geometrically connected components}

\begin{att}
Throughout this section, $K$ denotes a non-trivially valued field.
\end{att}

The results presented in this section are heavily inspired by the work of Jérôme Poineau in \cite{Poi14}. These results examine the variation of the connected components of a particular family of analytic spaces.

\subsection{Connectedness in Shilov's Spaces}

\begin{deff}[Shilov's space]
Let $X$ be a $K$-affinoid space, we say that $X$ is Shilov if $X$ is strict and if the Shilov boundary of $X$ reduces to a single point.
\end{deff}

\begin{nott}
Let $X$ be a $K$-affinoid Shilov space, we denote by $\PtShilov{X}$ the unique point on the Shilov boundary of $X$.
\end{nott}

\begin{lem}
Let $X$ be a $K$-affinoid Shilov space, then $X$ is connected.
\end{lem}
\begin{dem}
Since the Shilov boundary $X$ is the inverse image of the generic points of the irreducible components of the reduction of $X$, the space $\tilde{X}$ has a unique irreducible component, and $X$ is connected.
\end{dem}

\begin{prop}\label{PROP - Sur niveaux connexes sur espace de Shilov}
Let $X$ be a $K$-affinoid Shilov space and $f$ an analytic function on $X$., for every strictly positive real number $\varepsilon$, the set $\{ \varepsilon \leqslant \abs{f} \}$ is connected.
\end{prop}
\begin{dem}
According to Section 6.2.1/4 of \cite{BGR}, there exists a nonzero scalar $c$ in $K$ and a nonzero integer $m$ such that $cf^m$ has spectral norm equal to $1$. Even if we have to replace $f$ with $cf^m$, we can therefore assume that the spectral norm of $f$ is $1$, and the affinoid domains in question are then related by the following equation:
\[ \left\{ \varepsilon \leqslant \abs{cf^m} \right\} = \left\{ \left( {\varepsilon \over \abs{c}} \right)^{1 \over m} \leqslant \abs{f} \right\} \]
\begin{enumerate}
\item Suppose that $\varepsilon$ is strictly greater than $1$, the set $\{ \varepsilon \leqslant \abs{f} \}$ is empty, and therefore connected.
\item Suppose that $\varepsilon$ equals $1$, if $\{ 1 \leqslant \abs{f} \}$ is disconnected, then its reduction is also disconnected. However, the reduction of $\{ 1 \leqslant \abs{f} \}$ is the open set $\OuvZa{\tilde{f}}$, and since $\OuvZa{\tilde{f}}$ is an open set in the irreducible space $\tilde{X}$, it is connected. It follows that $\{ 1 \leqslant \abs{f} \}$ is connected.
\item Suppose that $\varepsilon$ is strictly less than $1$ and that $\{ \varepsilon \leqslant \abs{f} \}$ is not connected. Let $V$ be the connected component of $\{ \varepsilon \leqslant \abs{f} \}$ containing $\PtShilov{X}$ and $V_1$ the finite union of the other connected components of $\{ \varepsilon \leqslant \abs{f} \}$. Since $V_1$ is an affinoid domain of $X$ and does not contain the Shilov boundary of $X$, by the maximum principle given by Property 2.5.20 of \cite{Ber90}, the Shilov boundary of $V_1$ is included in the relative boundary of $V_1$ in $X$. Furthermore, since the relative boundary of $V_1$ in $X$ coincides with the topological boundary of $V_1$ in $X$ by Property 3.1.3 of \cite{Ber90}, the Shilov boundary of $V_1$ is contained in $\{ \abs{f} = \varepsilon \}$ and $\abs{f}$ is a constant equal to $\varepsilon$ on $V_1$. Since $X$ is compact, there exist two open sets $U$ and $U_1$ containing $V$ and $V_1$, respectively, such that $U \cap U_1$ is empty. Since the set $X \priv (U \cup U_1)$ is closed and nonempty by the connectedness of $X$, it is compact, and $\abs{f}$ attains a maximum $M$ on it. Let $\delta$ be a strictly positive real number such that $\varepsilon - \delta$ is strictly greater than $M$. The sets $V$ and $V_1$ are contained in $\{ \varepsilon - \delta \leqslant \abs{f} \}$, which is itself contained in $U \cup U_1$. Let $W_1$ denote the set $\{ \varepsilon - \delta \leqslant \abs{f} \} \cap U_1$. Since $W_1$ is affinoid, the same reasoning used earlier for $V_1$ shows that $\abs{f}$ is constant and equal to $\varepsilon - \delta$ on $W_1$. However, since $W_1$ contains $V_1$, this is a contradiction, and $\{ \varepsilon \leqslant \abs{f} \}$ is connected.
\end{enumerate}
\end{dem}

\begin{deff}[Geometric Shilov's space]
Let $X$ be a $K$-affinoid space, we say that $X$ is geometrically Shilov if, for every complete valued extension $L$ of $K$, the space $X_L$ is an $L$-affinoid Shilov space.
\end{deff}

\begin{prop}
Let $X$ be a $K$-affinoid space that is geometrically Shilov, and let $f$ be an analytic function on $X$. For every positive real number $\varepsilon$, the set $\{ \varepsilon \leqslant \abs{f} \}$ is geometrically connected.
\end{prop}
\begin{dem}
The result follows from the property \ref{PROP - Sur niveaux connexes sur espace de Shilov}.
\end{dem}

\subsection{Geometric connectedness in Shilov spaces}

\begin{deff}[Relative Shilov's space]
Let $X$ be a $K$-affinoid space and $Y$ a $K$-affinoid $X$-space, we say that $Y$ is Shilov on $X$ if:
\begin{enumerate}[label=$(SH_{\arabic*})$]
\item For every $x$ in $X$, the space $Y_x$ is a $\CrpResC{x}$-affinoid Shilov space.
\item There exists a universally multiplicative affinoid $X$-space $Z$ such that $Y$ is equal to $X \times_K Z$.
\end{enumerate}
\end{deff}

\begin{nott}
Let $X$ be a $K$-analytic space, $Y$ a $K$-analytic Shilov $X$-space, and $x$ a point in $X$. We denote by $\PtShilov{Y,x}$ the unique point on the Shilov boundary of $Y_x$.
\end{nott}

\begin{theo}\label{THEO - controle des composantes connexes relative}
Let $(X,x)$ be a pointed $K$-affinoid space, $Y$ a $K$-affinoid Shilov space over $X$, $y$ a point in $Y_x$, and $f$ an analytic function on $Y$ that is nonzero at $y$. For every nonzero real number $\varepsilon$ strictly less than $\abs{f(y)} $, there exists a neighborhood $U$ of $x$ such that $\{ \varepsilon \leqslant \abs{f} \} \times_X U$ has nonempty connected fibers over $U$ and contains $y$.
\end{theo}
\begin{dem}
Let $t$ be a point in $X$. Let $g_t$ denote the analytic function $f(t)$ on $Y_t$, $\delta_t$ the real number $\abs{g_t(\PtShilov{Y,t})}$, and $U$ the set of points $t$ in $X$ such that $\varepsilon$ is strictly less than $\delta_t$; by assumption, the set $U$ contains the point $x$. Furthermore, since the Shilov section of $X$ is continuous by (3.6) of \cite{Poi13}, the set $U$ is open, and for every $t$ in $U$, the fiber of $\{ \varepsilon \leqslant \abs{f} \} \times_X U$ at $t$, which is equal to $\{ \varepsilon \leqslant \abs{g_t} \}$, is connected by property \ref{PROP - Sur niveaux connexes sur espace de Shilov}. Finally, the fiber of $\{ \varepsilon \leqslant \abs{f} \} \times_X U$ at $x$ contains $y$ by assumption.
\end{dem}

\begin{deff}[Relative geometric Shilov's space]
Let $X$ be a $K$-affinoid space and $Y$ a $K$-affinoid $X$-space. We say that $Y$ is geometrically Shilov on $X$ if, for every $x$ in $X$, the space $Y_x$ is a $\CrpResC{x}$-affinoid space that is geometrically Shilov.
\end{deff}

\begin{prop}\label{PROP - controle des composantes geo connexes relative}
Let $(X,x)$ be a $K$-affinoid space with a base point, $Y$ a $K$-affinoid $X$-space that is geometrically Shilov, $y$ a point in $Y_x$, and $f$ an analytic function on $Y$ that is nonzero at $y$. For every nonzero real number $\varepsilon$ strictly less than $\abs{f(y)} $, there exists a neighborhood $U$ of $x$ such that $\{ \varepsilon \leqslant \abs{f} \} \times_X U$ has nonempty, geometrically connected fibers over $U$ and contains $y$.
\end{prop}
\begin{dem}
The result follows from the theorem \ref{THEO - controle des composantes connexes relative}.
\end{dem}

\subsection{Cutting connected components of fibers}

\begin{deff}[Relative cutting]
Let $X$ and $Y$ be two $K$-analytic spaces, and let $f$ be a morphism from $Y$ to $X$. An $X$-cutting of $Y$ is a $G$-covering of $Y$ by compact analytic domains whose fibers over $X$ are geometrically connected.
\end{deff}

\begin{theo}\label{THEO - Decoupage de composantes connexes de fibres}
Let $(X,x)$ be a strict pointed compact $K$-analytic space, $(Y,y)$ a $K$-analytic pointed $X$-space that is geometrically Shilov and quasi-smooth over $X$, $Z$ a Zariski closed set of $Y$ containing $y$ such that, for every $t$ in $X$, the set $Z_t$ has an empty interior in $Y_t$, and $V$ a compact analytic domain of $Y$ containing $y$, there exists an analytic domain $W$ in $Y$ and a neighborhood $U$ of $x$ such that $W \cap Z$ equals $V \cap Z$ and such that $W \times_X U$ has geometrically connected fibers over $U$ and contains $y$.
\end{theo}
\begin{dem}
By Theorem (3.11) of \cite{Duc20}, there exists a strictly affinoid, quasi-étale $X$-analytic space $\tilde{X}$ that is surjective onto $X$ such that $V \times_X \tilde{X}$ has a $\tilde{X}$-decomposition $(V_i)_{i \in I}$ into strictly analytic domains.
\[ \xymatrix{
&& \bigsqcup_{i \in I} V_i \ar@{>>}[d] \\
V \ar[d] && \ar@{>>}[ll] V \times_X \tilde{X} \cartesien \ar[dd]\\
Y \ar[d] \\
X && \ar@{>>}[ll] \tilde{X} } \]
Let $h$ be a nonzero analytic function that vanishes on $Z$. By restricting $X$ and by the continuity of the Shilov section, for every $t$ in $X$, the element $h_t(\PtShilov{Y,t})$ is nonzero, and the set $\{ \abs{h} = 0 \}$ is a nonempty interior set in each fiber. Let $i$ be an element of $I$. Since the analytic $X$-space $V_i$ is compact and quasi-smooth, Corollary (3.8) of \cite{Duc20} implies the existence of a compact analytic domain $W_i$ in $V_i$ that intersects every fiber of $V_i$ and is disjoint from $\{ \abs{h} = 0 \} \times_X V_i$. Let $W$ denote the union of the images in $Y$ of the sets in the family $(W_i)_{i \in I}$; the set $W$ is a strictly analytic compact domain, does not intersect $\{ \abs{h} = 0 \}$, and intersects all connected components of the fibers of $V$ over $X$.
\[ \xymatrix{
W \ar[dd] && \bigsqcup_{i \in I} W_i \ar@{>>}[ll] \ar[d] \\
&& \bigsqcup_{i \in I} V_i \ar@{>>}[d] \\
V \ar[d] && \ar@{>>}[ll] V \times_X \tilde{X} \cartesien \ar[dd] \\
Y \ar[d]_f \\
X && \ar@{>>}[ll] \tilde{X} } \]
Since $h$ is a nonzero continuous function on the compact set $W$, the function $\abs{h}$ is bounded on $W$ and attains its nonzero minimum $\varepsilon$ there. Now, since the set $W$ contains all the connected components of the fibers of $V$ over $X$ and is contained in $\{ \varepsilon \leqslant \abs{h} \} $, there exists a neighborhood $U$ of $x$ such that the fibers of $(V \cup \{ \varepsilon \leqslant \abs{h} \}) \times_X U$ are connected and therefore geometrically connected by the property \ref{PROP - controle des composantes geo connexes relative}.
\end{dem}

\newpage \section{Raynaud-Gruson's dévissage}

\begin{att}
Throughout this section, $K$ denotes a non-trivially valued field.
\end{att}

\begin{deff}[Raynaud-Gruson's dévissage]
Let $(X,x)$ and $(Y,y)$ be two strict $K$-analytic pointed spaces, $f$ a pointed morphism from $(Y,y)$ to $(X,x)$, $\cF$ a coherent $\cOf{Y}$-module such that $y$ belongs to $\SuppFa{\cF}$, $r$ a positive nonzero integer, and $(n_i)_{1 \leqslant i \leqslant r}$ a strictly decreasing family of natural numbers, we call the following data a Raynaud-Gruson's $X$-dévissage of $\cF$ at $y$ in dimensions $(n_i)_{1 \leqslant i \leqslant r}$:
\begin{enumerate}
\item An analytic neighborhood $V$ of $y$ in $Y$.
\item For every $i$ between $1$ and $r$, a strict $K-$affinoid domain $T_i$ of a smooth $X-$space of pure relative dimension $n_i$ such that $T_i$ has geometrically connected fibers on $X$.
\item For every $i$ between $1$ and $r$, a point $t_i$ in $T_i$ above $x$.
\item For every $i$ between $1$ and $r$, a free coherent $\cOf{T_i}$-module $\cL_i$.
\item For every $i$ between $1$ and $r-1$, a coherent $\cOf{T_i}$-module $\cP_i$ such that $t_i$ belongs to the support of $\cP_i$.
\item A coherent $\cOf{T_r}$-module $\cP_r$ such that, if $\cP_r$ is nonzero, the point $t_r$ belongs to the support of $\cP_r$.
\item A finite $X$-morphism $\pi_1$ from the support of $\cF_V$ to $T_1$ such that $\pi_1(y) = t_1$.
\item For every $i$ between $2$ and $r$, a finite $X$-morphism $\pi_i$ from the support of $\cP_{i-1}$ to $T_i$ such that $\pi_1(t_{i-1}) = t_i$.
\item A morphism $\alpha_1$ from $\cL_1$ to $(\pi_1)_\ast(\cF_V)$ whose conkernel is $\cP_1$ and such that:
\[ \adh[(T_1)_x]{ \BijL{\Tige{\alpha_1}{x}} } = (T_1)_x \] 
\item For every $i$ between $2$ and $r$, a morphism $\alpha_i$ from $\cL_i$ to $(\pi_i)_\ast(\cP_{i-1})$ whose cokernel is $\cP_i$ and such that:
\[ \adh[(T_i)_x]{ \BijL{\Tige{\alpha_i}{x}} } = (T_i)_x \]
\end{enumerate}
We then denote $(V,\{T_i,t_i,\pi_i,\alpha_i,\cL_i,\cP_i\}_{1 \leqslant i \leqslant r})$ these data for the Raynaud-Gruson's $X$-dévissage of $\cF$ at $y$ in dimensions $(n_i)_{1 \leqslant i \leqslant r}$.
\end{deff}

\begin{voc}[Total Raynaud-Gruson's dévissage]
Let $(X,x)$ and $(Y,y)$ be two pointed strict $K$-analytic spaces, $f$ a pointed morphism from $(Y,y)$ to $(X,x)$, $\cF$ a coherent $\cOf{Y}$-module such that $y$ belongs to $\SuppFa{\cF}$, $r$ a positive nonzero integer, $(n_i)_{1 \leqslant i \leqslant r}$ a strictly decreasing family of natural numbers, and $\cD = (V,\{T_i,t_i,\pi_i,\alpha_i,\cL_i,\cP_i\}_{1 \leqslant i \leqslant r})$ a Raynaud-Gruson's $X$-dévissage of $\cF$ at $y$ in dimension $(n_i)_{1 \leqslant i \leqslant r}$, we say that $\cD$ is total if $\cP_r$ is zero.
\end{voc}

Let $(X,x)$ and $(Y,y)$ be two pointed strict $K$-analytic spaces, $f$ a pointed morphism from $(Y,y)$ to $(X,x)$, $\cF$ a coherent $\cOf{Y}$-module such that $y$ belongs to $\SuppFa{\cF}$, $r$ a positive nonzero integer, $(n_i)_{1 \leqslant i \leqslant r}$ a strictly decreasing family of natural numbers, and $(V,\{T_i,t_i,\pi_i,\alpha_i,\cL_i,\cP_i\}_{1 \leqslant i \leqslant r})$ a Raynaud-Gruson's $X$-dévissage of $\cF$ at $y$ in dimension $(n_i)_{1 \leqslant i \leqslant r}$, we have the following diagram and exact sequences:
{\small~ \[ \xymatrix{
\cL_1 \to (\pi_1)_\ast(\cF_V) \to \cP_1 \to 0 &&& (\SuppFa{\cF_V},y) \ar[d]^{\pi_1} \ar@{^{(}->}[r] & (V,y) \ar[dd] & \\
\cL_2 \to (\pi_2)_\ast(\cP_1) \to \cP_2 \to 0 \ar@{.}[dd] && (\SuppFa{\cP_1},t_1) \ar@{^{(}->}[r] \ar[d]^{\pi_2} & (T_1,t_1) \ar[rddd] &&\\
&& (T_2,t_2) \ar[rrdd] && (Y,y) \ar[dd] \\
\cL_{r-1} \to (\pi_{r-1})_\ast(\cP_{r-2}) \to \cP_{r-1} \to 0 & (\SuppFa{\cP_{r-1}},t_{r-1}) \ar[d]^{\pi_r} \ar@{^{(}->}[r] & (T_{r-1},t_{r-1}) \ar[rrd] \ar@{.}[u] &&& \\
\cL_r \to (\pi_r)_\ast(\cP_{r-1}) \to \cP_r = 0 & (T_r,t_r) \ar[rrr] &&& (X,x)} \] }

Unlike the dévissages presented in \cite{Duc18} and \cite{Duc20}, the notion of a Raynaud-Gruson's dévissage given here may seem weaker. We require only an analytic neighborhood of the point $y$ in $Y$, and not an affinoid neighborhood. Similarly, for any $i$ between $1$ and $r$, we impose no conditions on the fiber of $\pi_i$ at the point $t_i$. This generality is, in reality, only apparent, since we require the $X$-spaces $(T_i)_{1 \leqslant i \leqslant r}$ to be quasi-smooth.

\begin{theo}\label{THEO - existence de devissage connexe}
Let $(X,x)$ and $(Y,y)$ be two good pointed strict $K$-analytic spaces, $f$ a pointed morphism from $(Y,y)$ to $(X,x)$ such that $y$ lies in the interior of $f$, $\cF$ a coherent $\cOf{Y}$-module such that $y$ belongs to $\SuppFa{\cF}$, and $n$ the integer $\dimAn[K]{y}{\Tige{\cF}{x}}$. There exists a Raynaud-Gruson's $X$-dévissage of $\cF$ at $y$ in dimensions less than $n$.
\end{theo}
\begin{dem}
Since $K$ is not trivially valued and the result is local on $X$ and $Y$, we may assume that $X$ and $Y$ are strict $K$-affinoid. Furthermore, since $\SuppFa{\cF}$ is Zariski-closed in $Y$, the interior of $f$ is contained in the interior of $\SuppFa{\cF}$ by (1.5.5) of \cite{Ber93} and the fact that a closed immersion is boundaryless. By the property \ref{PROP - factorisation fini et lisse connexe sur un ferme}, there exists a strict affinoid rational neighborhood $V$ of $y$ in $Y$, an affinoid neighborhood $W$ of $x$ in $X$, and an $n$-polyradius $r$ in $\abs{K^\times}\! {}^{\eQ}$ such that $\rest{f}{(V \cap \SuppFa{\cF})}$ factors through a finite morphism $\pi$ from $V \cap \SuppFa{\cF}$ to $\eD^n_W(0,r)$. Let $Z$ denote the affinoid space $V \cap \SuppFa{\cF}$, $T$ the polydisk $\eD^n_W(0,r)$, and $t$ the image of $y$ under $\pi$.
\[ \xymatrix{
(\SuppFa{\cF},y) \ar[d] && \ar[ll] (Z,y) \ar[d]^{\pi} \\
(Y,y) \ar[d] && (T,t) \ar[d] \\
(X,x) && \ar[ll] (W,x)
} \]
Let $\cG$ denote the coherent $\cOf{T}$-module $\pi_\ast(\cF_Z)$ and $\delta$ the Krull dimension of $\cOfg{Z_x}{y}$. By the finiteness of $Z_x$ over $T_x$ and by (3.2.4) in \cite{Duc18}, we have:
\[ \dimCtr{T_x,t} = \dimCtr{Z_x,y} \]
However, as:
\[ \dimAn[K]{y}{Z_x} = \dimAn[t]{K}{T_x} = n \]
we have :
\begin{equation}
\delta = \dimKr{\cOfg{T_x}{t}} \label{EQ - egalite de dim de krull}
\end{equation}
by (3.2.9) of \cite{Duc18}. Furthermore, since the support of $\cG$ is equal to $\pi\fens{Z}$ and since $Z_x$ has analytic dimension $n$ in $y$, we have:
\begin{equation}
\dimAn[K]{t}{\Tige{\cG}{x}} = \dimAn[K]{t}{\pi\fens{Z_x}} = n \label{EQ - Dim locale sur fibre}
\end{equation}
by (1.5.10) of \cite{Duc18}. Since the fibre $T_x$ is a $\CrpResC{x}$-analytic space of pure dimension $n$, the equality (\ref{EQ - Dim locale sur fibre}) implies that $\SuppFa{\cG}_x$ contains the connected component of $t$ in $T_x$ and that the annihilator of $\Tige{\cG}{x}$ is zero at $t$. It follows that:
\begin{equation}
\delta = \dimKr{\Tige{\cG}{x}_t} \label{EQ - dim krull locale sur fibre de module}
\end{equation}
\noindent Let us reason by recurrence with respect to $n$ and denote by $c$ the integer $\Coprofl[t]{\Tige{\cG}{x}}$.
\begin{enumerate}[itemindent=*,leftmargin=*]
\item If $n$ is zero, then the integer $c$ is zero by (8.2.2) of \cite{Duc18}, and thus $\Tige{\cG}{x}_t$ is a $\cOfg{T_x}{t}$-module of finite type with zero codepth and maximal Krull dimension over a regular local ring, and is therefore free by ($0_{IV}$, \S 17.3.4) of \cite{EGA4}. Let $(f_i)_{1 \leqslant i \leqslant r}$ be a family of global sections of $\cG$ whose image in $\Tige{\cG}{t}$ is a basis, $\cL$ the $\cOf{T}$-free module $\cOf[r]{T}$ and $\alpha$ the morphism from $\cL$ to $\cG$ defined by $(f_i)_{1 \leqslant i \leqslant r}$; by Nakayama’s lemma, the morphism $\alpha$ is surjective in $t$. Furthermore, the restriction of $\alpha$ to $T_x$ is bijective in $t$ because it is a surjective morphism between two $\cOfg{T_x}{t}$-modules of finite types of the same rank. Let $F_\alpha$ be the Zariski closure of the non-bijectivity of $\Tige{\alpha}{x}$ in $T_x$, $f$ an analytic function on $T_x$ vanishing in $F_\alpha$ such that $f(t)$ is non-zero, and $\eta$ a positive, non-zero real number strictly less than $\abs{f(t)} $, even if we have to restrict $X$ around $x$, there exists an analytic function $g$ on $T$ such that its restriction $g_0$ to $T_x$ satisfies:
\[ \norm{ g_0 - f } \leqslant \eta \]
Let $\varepsilon$ be a real number strictly less than $\abs{f(t)}$ and strictly greater than $\eta$, and let $W_\epsilon$ denote the set:
\[ W_\varepsilon = \{ \abs{g} \geqslant \varepsilon \} \]
Since $T$ is geometrically Shilov on $X$ and since the desired result is local on $X$, even if this means restricting $X$ to the neighbourhood of $x$, the space $W_\varepsilon$ contains $t$, has geometrically connected fibres on $X$ by property \ref{PROP - controle des composantes geo connexes relative}, and does not contain $F_\alpha$. Provided we restrict to $W_\varepsilon$, we may assume that $\Tige{\alpha}{x}$ is bijective. Furthermore, since the support of the kernel of $\alpha$ is closed in $Y$, provided that $X$ is restricted to the neighbourhood of $x$, we may assume that $\alpha$ is bijective. Ve then obtain a Raynaud-Gruson's $X$-dévissage of $\cF$ at $y$ in dimension $n$.
\item If $n$ is non-zero, let us assume that the result has been proved for dimensions strictly less than $n$.
\begin{att}
In the remainder of the proof, we shall assign the index $x$ to those points of $T$ that are regarded as elements of $T_x$.
\end{att}
\begin{itemize}
\item Let $z_x$ be a point in $T_x$ such that $\Algbr{z}_x$ is the generic point of the connected component of $\Algbr{(T_x)}$ containing $\Algbr{t}_x$. Since the support of $\Tige{\cG}{x}$ contains the connected component of $t$, the $\CrpRes{\Algbr{z}_x}$-vector space $\Tige{(\Tige{\Algbr{\cG}}{\Algbr{x}})}{\Algbr{z}_x}$ is of non-zero dimension $r$. Let $(f_i)_{1 \leqslant i \leqslant r}$ be a family of global sections of $\cG$ whose image in $\Tige{(\Tige{\Algbr{\cG}}{\Algbr{x}})}{\Algbr{z}_x}$ forms a basis, $\cL$ the free $\cOf{T}$-module $\cOf[r]{T}$ and $\alpha$ the morphism from $\cL$ to $\cG$ defined by $(f_i)_{1 \leqslant i \leqslant r}$, as we have:
\[ \cOfg{\Algbr{(T_x)}}{\Algbr{z}_x} = \CrpRes{\Algbr{z}_x} \]
The morphism $\Tige{\Algbr{\alpha}}{\Algbr{x}}_{\Algbr{z}_x}$ is an isomorphism and $z$ belongs to $\BijL{\Tige{\alpha}{x}}$. Now, since $t$ is a Zariski specialisation of $z$, it belongs to $\adh[T_x]{ \BijL{\Tige{\alpha}{x}} }$.
\item Since the scheme $\Algbr{(T_x)}$ is regular, the local ring of $\Algbr{(T_x)}$ at $\Algbr{t}_x$ is an integral ring whose field of fractions is $\CrpRes{\Algbr{z}_x}$. Now, since $\cL$ is a free $\cOf{T}$-module, the morphism from $\Tige{\Algbr{\cL}}{\Algbr{x}}_{\Algbr{t}_x}$ to $\Tige{\Algbr{\cL}}{\Algbr{x}}_{\Algbr{z}_x}$ is injective. From this, by the commutative diagram, we deduce:
\[ \xymatrix{ \Tige{\Algbr{\cL}}{\Algbr{x}}_{\Algbr{t}_x} \ar[rr]^{\Tige{\Algbr{\alpha}}{\Algbr{x}}_{\Algbr{t}_x}} \ar[d] && \Tige{\Algbr{\cG}}{\Algbr{x}}_{\Algbr{t}_x} \ar[d] \\
\Tige{\Algbr{\cL}}{\Algbr{x}}_{\Algbr{z}_x} \ar[rr]^{ \Tige{\Algbr{\alpha}}{\Algbr{x}}_{\Algbr{z}_x} } && \Tige{\Algbr{\cG}}{\Algbr{x}}_{\Algbr{z}_x} } \]
that the morphism $\Tige{\Algbr{\alpha}}{\Algbr{x}}_{\Algbr{t}_x}$ is injective and therefore that $\Tige{\alpha}{x}_t$ is injective.
\begin{enumerate}
\item Suppose that $c$ is zero, the morphism $\Tige{\alpha}{x}_t$ is surjective and $t$ belongs to $\BijL{\Tige{\alpha}{x}}$. Let $f$ be an analytic function on $T_x$ vanishing on the support of $\cP_x$ such that $f(t)$ is non-zero, and let $\eta$ be a positive, non-zero real number strictly less than $\abs{f(t)}$, even if this requires restricting $X$ to a neighbourhood of $x$, there exists an analytic function $g$ on $T$ such that its restriction $g_0$ to $T_x$ satisfies:
\[ \norm{ g_0 - f } \leqslant \eta \]
Let $\varepsilon$ be a real number strictly less than $\abs{f(t)}$ and strictly greater than $\eta$, and let $W_\epsilon$ denote the set:
\[ W_\varepsilon = \{ \abs{g} \geqslant \varepsilon \} \]
Since $T$ is geometrically Shilov on $X$ and since the desired result is local on $X$, even if this means restricting $X$ to the neighbourhood of $x$, the space $W_\varepsilon$ contains $t$, has geometrically connected fibres on $X$ by property \ref{PROP - controle des composantes geo connexes relative} and does not contain the support of $\cP_x$ in $T_x$. Even if we restrict ourselves to $W_\varepsilon$, we may assume that $\Tige{\alpha}{x}$ is bijective. Furthermore, since the support of the kernel of $\alpha$ is closed in $Y$, provided we restrict $X$ to the neighbourhood of $x$, we may assume that $\alpha$ is bijective, and we then obtain a Raynaud-Gruson's $X$-dévissage of $\cF$ at $y$ in dimensions $n$.
\item Let us assume that $c$ is non-zero. 
\begin{itemize}
\item The morphism $\Tige{\alpha}{x}_t$ is not surjective because, otherwise, since the local ring of $\Algbr{(T_x)}$ at $\Algbr{t}_x$ is regular, the codepth of $\Tige{\cG}{x}_t$ would be zero, which is absurd. Let $\cP$ denote the cokernel of $\alpha$. Since $\Tige{\alpha}{x}_t$ is not surjective, the point $t$ lies in the support of $\cP_x$. Now, since the support of $\cP_x$ is contained in $T_x \priv \BijL{\Tige{\alpha}{x}}$, we have:
\[ \dimAn[K]{t}{ \SuppFa{\Tige{\cP}{x}} } \leqslant n-1 \]
This ensures that $\cP$ has a Raynaud-Gruson's $X-$dévissage :
\[ \cD = (V_\cP,\{T_i,t_i,\pi_i,\alpha_i,\cL_i,\cP_i\}_{1 \leqslant i \leqslant r}) \]
in $t$ in dimensions lower than $n$.
\item Since the set $V_\cP$ is an analytic neighbourhood containing $t$, by (4.9) in \cite{Duc07}, even if we restrict $X$ to the neighbourhood of $x$, the fibres of $V_\cP$ have empty interiors and there exists an analytic domain $W$ such that $V \cap \SuppFa{\cP}$ equals $W \cap \SuppFa{\cP}$ and such that $W$ has geometrically connected fibres on $X$ by the theorem \ref{THEO - Decoupage de composantes connexes de fibres}. Let $V_W$ denote the set $\reci{\pi}\fens{W} \cap V$ and $\pi_W$ the restriction of $\pi$ to $\reci{\pi}\fens{W}$, since $W$ contains $t$, the data of $\cD$ together with the data $(V_W,\{W,t,\pi_W,\rest{\alpha}{W},\rest{\cL}{W},\rest{\cP}{W}\})$ then define a Raynaud-Gruson's $X$-dévissage of $\cF$ at $y$ in dimensions less than $n$.
\end{itemize}
\end{enumerate}
\end{itemize}
\end{enumerate}
\end{dem}

\newpage \section{Enjoliveurs}

\begin{att}
Throughout this section, $K$ denotes a non-trivially valued field, $\Gamma$ a subgroup of $\eR_{\PM +}^\ast$ such that $\Gamma\absNP{K}$ is not the trivial group, $\fE$ a concrete category over the category of $\Gamma$-strict $K$-analytic spaces and $\fC$ a class of morphisms in $\fE$ containing the isomorphisms. Hereafter, whenever $\Gamma = \eR_{\PM +}^\times$, we shall omit this notation.
\end{att}

\begin{nott}
Let $W$ be an element of $\fE$, we denote by $\Conc{W}$ the image of $W$ in the category of $\Gamma$-strict $K$-analytic spaces.
\end{nott}

\subsection{Admissible presentation and enjoliveurs}

\begin{nott}
We define, recursively for a natural number $h$, a class $\DiagApplAn[h]{K}{\Gamma}(\fC,\fE)$ of pairs consisting of a $K$-analytic $\Gamma$-strict space and a morphism from $\fE$ by the following procedure:
\begin{enumerate}
\item Let $\DiagApplAn[0]{K}{\Gamma}(\fC,\fE)$ denote the class of pairs:
\[ (S,\xymatrix{ Y \ar[rr]^v && X }) \]
such that:
\begin{itemize}
\item The morphism $v$ is an isomorphism in $\fE$.
\item The $K$-analytic space $S$ is the divisor $\emptyset$ of $\Conc{Y}$.
\end{itemize}
\item Let $h$ be a non-zero integer, the pair:
\[ (S,\xymatrix{ Y \ar[rr]^v && X }) \]
belongs to $\DiagApplAn[h]{K}{\Gamma}(\fC,\fE)$ if and only if there exists:
\begin{enumerate}
\item A couple:
\[ (T,\xymatrix{ Z \ar[rr]^u && X }) \]
belonging to $\DiagApplAn[h-1]{K}{\Gamma}(\fC,\fE)$.
\item An object $W$ of $\fE$ over $Z$ whose structural morphism belongs to $\fC$.
\item A closed $K$-analytic subspace $\Sigma$ of $\Conc{W}$ factoring through $\Conc{W} \times_{\Conc{Z}} T$, such that $v$ has a factorisation:
\[  \xymatrix{ v : Y \ar[rr] && W \ar[rr] && Z \ar[rr]^u && X } \]
identifying $\Conc{Y}$ to $\BLup{\Sigma}{\Conc{W}}$ and $S$ to $\BLupD{\Sigma}{\Conc{W}}$.
\end{enumerate}
\item Let $\DiagApplAn[\infty]{K}{\Gamma}(\fC,\fE)$ denote the union of $(\DiagApplAn[h]{K}{\Gamma}(\fC,\fE))_{h \in \eN}$.
\end{enumerate}
\end{nott}

\begin{deff}[Admissible presentation]
Let $X$ and $Z$ be two objects in $\fE$ and $S$ a divisor of $\Conc{Z}$, an $(\fC,\fE)$-admissible $\Gamma$-presentation of $(Z,S)$ over $X$ is a pair:
\[ (S,\xymatrix{ Z \ar[rr]^v && X }) \]
belonging to $\DiagApplAn[\infty]{K}{\Gamma}(\fC,\fE)$.
\end{deff}

Let $X$ and $Z$ be two objects in $\fE$ and $S$ a divisor of $\Conc{Z}$, an $(\fC,\fE)$-admissible $\Gamma$-presentation of $(Z,S)$ is given by a morphism $\varphi$ from $Z$ to $X$ in $\fE$ that has a factorization:
\[ \varphi : \xymatrix{ Z = Z_{2n} \ar[rr] && Z_{2n-1} \ar[rr] && \ldots \ar[rr] && Z_0 = X } \]
and, for every $i$ between $0$ and $2n$, a closed $K$-analytic subspace $S_i$ of $\Conc{Z_i}$ such that:
\begin{enumerate}
\item We have $S_{2n} = S$ and $S_0 = \emptyset$.
\item If $i$ is even, the morphism from $Z_{i+1}$ to $Z_i$ belongs to $\fC$, and $S_{i+1}$ factors through $\Conc{Z_{i+1}} \times_{\Conc{X}} S_i$.
\item If $i$ is odd, the space $\Conc{Z_{i+1}}$ is $\BLup{S_i}{\Conc{Z_i}}$, and $S_{i+1}$ is the exceptional divisor $\Conc{Z_{i+1}} \times_{\Conc{X}} S_i$.
\end{enumerate}

\begin{deff}[Enjoliveur]
Let $X$ be an object in $\fE$, an $(\fC,\Gamma)$-enjoliveur of $X$ is the specification of a pair $(Z,S)$, where $Z$ is an object in $\fE$ and $S$ is a divisor of $\Conc{Z}$, which has an $(\fC,\fE)$-presentation that is $\Gamma$-admissible.
\end{deff}

\begin{voc}[Enjoliveur's class]
Let $X$ be an object in $\fE$ and $(Z,S)$ a $(\fC,\Gamma)$-embellisher of $X$, we call class of $(Z,S)$ the smallest integer $n$ such that $(Z,S)$ has an $(\fC,\fE)$-admissible $\Gamma$-presentation belonging to $\DiagApplAn[n]{K}{\Gamma}(\fC,\fE)$.
\end{voc}

\subsection{Flattening with étale enjoliveurs}

\begin{nott}
Let $X$ and $Y$ be two $K$-analytic spaces, $f$ a morphism from $Y$ to $X$ and $\cF$ a coherent $\cOf{Y}$-module.
\begin{enumerate}
\item For any integer $n$, we denote by $\LNPlat[\geqslant n]{\cF}{X}$ the set of points in $\LNPlat{\cF}{X}$ at which the latter has relative dimension greater than $n$ over $X$.
\item For any integer $n$, we denote by $\LNPlatS[\geqslant n]{\cF}{X}$ the saturated set of $\LNPlat[\geqslant n]{\cF}{X}$ with respect to the equivalence relation defined by $f$.
\end{enumerate}
\end{nott}

\begin{deff}[Flat module in controled dimensions]
Let $X$ and $Y$ be two $K$-analytic spaces, $f$ a morphism from $Y$ to $X$ and $\cF$ a coherent $\cOf{Y}$-module, we say that $\cF$ is $X$-flat in dimensions greater than $n$ if $\LNPlat[\geqslant n]{\cF}{X}$ is empty
\end{deff}

Theorem (6.6), obtained by Antoine Ducros in \cite{Duc20}, although it deals only with morphisms between compact spaces, is not exactly the same statement -mutatis mutandis- as Theorem (5.2.2) of \cite{RG71}, which was obtained in scheme theory. Due to the wild behavior of the Zariski topology in analytic geometry, the flattening procedure must, at a minimum, allow for immersions of compact analytic domains in addition to blow-ups. More generally, one must resort to quasi-étale morphisms with compact sources, which are likely unavoidable according to (7.12) in \cite{Duc20}.

\begin{deff}[Étale enjoliveur]
Let $X$ be an object in $\fE$ and $(Z,S)$ a $(\fC,\Gamma)$-enjoliveur of $X$, we say that $(Z,S)$ is a étale $\Gamma$-enjoliveur of $X$ if $\fE$ is the category of compact $K$-analytic spaces and if $\fC$ is the class of quasi-étale morphisms.
\end{deff}

\begin{theo}[Theorem 6.6 of \cite{Duc20}]\label{THEO - Applatissement faible}
Let $X$ and $Y$ be two compact $\Gamma$-strict $K$-analytic spaces, $N$ and $n$ two integers such that $n$ is less than $N$, $f$ a morphism from $Y$ to $X$ and $\cF$ a coherent $\cOf{Y}$-module $X$-flat in dimensions greater than $N$, there exists an étale $\Gamma$-enjoliveur $(Z,S)$ of $X$ of class less than $N - n$ and a compact $\Gamma$-strict analytic domain $V$ of $Y \times_X Z$ such that the following conditions are satisfied:
\begin{enumerate}
\item The module $\cF_V \oslash_Z S$ is $Z$-flat in dimensions greater than $n$.
\item The set $V \priv (V \times_Z S)$ is contained in the inverse image of $Y \priv \LNPlat[\geqslant n]{\cF}{X}$.
\item The set $V \times_Z S$ is contained in the inverse image of $\LNPlatS[\geqslant n]{\cF}{X}$.
\item The image of $V$ in $Y$ contains $Y \priv \LNPlatS[\geqslant n]{\cF}{X}$.
\item The image of $V \times_Z S$ in $Y$ contains every point in $\LNPlat[\geqslant n]{\cF}{X}$ whose image on $X$ is not adherent to:
\[ f\fens{ \LNPlat[\geqslant n]{\cF}{X} } \cap \Abhy{X} \]
\item The image of $S$ on $X$ is contained in that of $\SuppFa{\cF}$ and, if $\cF$ is nonzero, this image has dimension less than $\dimAn{K}{\SuppFa{\cF}} - n$.
\end{enumerate}
\end{theo}

For the reader’s convenience, before presenting our result, let us review a few special cases, already discussed in \cite{Duc20}, on which several results concerning the embellishment of morphisms will be based. Let $X$ and $Y$ be two $\Gamma$-strict $K$-analytic spaces, $N$ and $n$ two integers such that $n$ is less than $N$, $f$ a morphism from $Y$ to $X$ and $\cF$ a coherent $\cOf{Y}$-module $X$-flat in dimensions greater than $N$. 
\begin{enumerate}
\item \textit{Suppose that $\cF$ is already $X$-flat in dimensions higher than $n$}\\
The theorem is therefore trivial, since the family consisting of the single pair:
\[ ( \xymatrix{ \emptyset \ar[rr] && X \ar[rr]^{\Id{X}} && X } , Y) \]
satisfies by definition the conditions of the theorem \ref{THEO - Applatissement faible}.
\item \textit{Let us assume that $n$ is zero}\\
Let $(Z,S)$ be the $\Gamma-$enjoliveur of $X$ of class less than $N - n$ given by theorem \ref{THEO - Applatissement faible}, the result obtained is a result of total flattening. In this case, the coherent sheaf $\mathscr F_Y \oslash_Z S$ is $Z-$flat.
\item \textit{Suppose that $\LNPlatS[\geqslant n]{\cF}{X}$ has an empty interior in $Y$}\\
The fourth statement of theorem \ref{THEO - Applatissement faible}, combined with the compactness of $V$, then implies that the image of $V$ in $Y$ is equal to $Y$.
\item \textit{Let us assume that $\LNPlatS[\geqslant n]{\cF}{X}$ is equal to $Y$}\\
The fourth statement of theorem \ref{THEO - Applatissement faible} is therefore trivial and provides no information about the image of $V$ in $Y$, whilst the second statement ensures that $V$ is contained in the reciprocal image of $S$, which implies that $\cF_V \oslash_Z S$ is zero. The first statement is therefore also obvious. The theorem \ref{THEO - Applatissement faible} is thus, in this case, devoid of content as regards flatness. This is consistent with the general principle mentioned in the introduction, according to which flattening techniques are not intended to create flatness ex nihilo, but rather to propagate it slightly where it already exists. Note that the previous point is also an illustration of this principle, but in the opposite case where there is already "a great deal" of flatness to begin with.\\
However, theorem \ref{THEO - Applatissement faible} is not trivial in this case. As we have just seen, it asserts the existence of a set-theoretic factorisation of the morphism from $V$ to $Z$ by $S$, which is a Zariski closed set of $Z$ with empty interior. If, moreover $X$ is reduced, then the fibre of $f$ at a point in $\Abhy{X}$ is empty, so that the image of $f$ contains no points in $\Abhy{X}$. Combining the first assertion with the second then implies that the image of $V$ in $Y$ is equal to $Y$.
\end{enumerate}

\subsection{Local enjoliveurs}

\begin{att}
In this subsection, $\fE$ denotes the category of $K$-analytic spaces.
\end{att}

\begin{voc}[Domanial morphism]
Let $X$ and $Y$ be two $K$-analytic spaces and $f$ a morphism from $Y$ to $X$, we say that $f$ is domanial if $f$ is a finite disjoint union of inclusions of affinoid domains.
\end{voc}

\begin{deff}[Local enjoliveur]
Let $(X,x)$ be a pointed $K$-analytic space and $(Z,S)$ a $\fC$-embellisher of $X$, we say that $(Z,S)$ is a local enjoliveur of $(X,x)$ if $\fC$ is the class of domanial morphisms and if the first step of the admissible presentation of $(Z,S)$ is the inclusion of an affinoid neighbourhood of $x$ into $X$.
\end{deff}

Let $X$ and $Z$ be two $K$-analytic spaces, $x$ a point in $X$ and $S$ a divisor of $Z$, the data of a $(\fC,\fE)$-admissible presentation of $(Z,S)$, such that $(Z,S)$ is a local enjoliveur of $(X,x)$, is the data of a morphism $\varphi$ from $Z$ to $X$ which admits a factorisation:
\[ \varphi : \xymatrix{ Z = Z_{2n} \ar[rr] && Z_{2n-1} \ar[rr] && \ldots \ar[rr] && Z_0 = X } \]
and, for every $i$ between $0$ and $2n$, a closed $K$-analytic subspace $S_i$ of $Z_i$ such that:
\begin{enumerate}
\item We have $S_{2n} = S$ and $S_0 = \emptyset$.
\item The morphism from $Z_1$ to $Z$ is given by the inclusion of an affinoid neighbourhood of $x$ into $X$.
\item If $i$ is even, the morphism from $Z_{i+1}$ to $Z_i$ is domanial and $S_{i+1}$ factors through $Z_{i+1} \times_{\Conc{X}} S_i$.
\item If $i$ is odd, the space $\Conc{Z_{i+1}}$ is $\BLup{S_i}{\Conc{Z_i}}$, and $S_{i+1}$ is the exceptional divisor $\Conc{Z_{i+1}} \times_{\Conc{X}} S_i$.
\end{enumerate}

\subsection{Base change by a local enjoliveur}

Let $(X,x)$ and $(Y,y)$ be two pointed strict $K$-analytic spaces, $f$ a pointed morphism from $(Y,y)$ to $(X,x)$, $(Z,S)$ a local enjoliveur of $(X,x)$, $U$ an analytic neighbourhood of $y$ in $Y$ and $V$ the $K$-analytic space $U \times_X Z$, we are here interested in the properties of the $Z$-space $V \oslash_Z S$.
\begin{itemize}
\item Since the morphism from $Y \times_X Z$ to $Y$ is a composition of domanial morphisms and blow-ups, the space $V \priv V \times_Z S$ is reduced as soon as $Y$ is reduced. Consequently, the morphism from $V \oslash_Z S$ to $Y$ is generically domanial and $V \oslash_Z S$ is reduced as soon as $Y$ is reduced.
\item Let $v$ be a point in $\Abhy{V \oslash_Z S}$, the morphism from $V \oslash_Z S$ to $Y$ is a local isomorphism at $v$, so that the image $y$ of $v$ in $Y$ belongs to $\Abhy{Y}$ and satisfies the following equality:
\[ d_K(v) = d_K(y) \]
\end{itemize}

\begin{lem}\label{LEM - dimension et enjoliveurs locaux}
Let $(X,x)$ and $(Y,y)$ be two pointed strict $K$-analytic spaces, $f$ a pointed morphism from $(Y,y)$ to $(X,x)$, $(Z,S)$ a local enjoliveur of $(X,x)$, $U$ an analytic neighbourhood of $y$ in $Y$ and $V$ the $K$-analytic space $U \times_X Z$.
\begin{enumerate}
\item If $Y$ is purely $d$-dimensional, then $V \oslash_Z S$ is purely $d$-dimensional.
\item If $f$ is generically of dimension\footnote{That is, the dimension of $f$ is equal to $d$ at every point of $\Abhy{Y}$.} $d$, then the morphism from $V \oslash_Z S$ to $Z$ is generically of dimension $d$.
\end{enumerate}
\end{lem}
\begin{dem}
Let $v$ be a point in $\Abhy{V \oslash_Z S}$. Since $v$ lies above a point in $\Abhy{Y}$, the dimension of the morphism from $Y \times_X Z$ to $Z$ at $v$ is equal to $d$, and the same holds for the dimension of the morphism from $V$ to $Z$ at $v$. Since $v$ belongs to $\Abhy{V \oslash_Z S}$, it has an open neighborhood in $V$ that is also an open set in $V \oslash_Z S$, so the dimension of the morphism from $V \oslash_Z S$ to $Z$ at $v$ is again equal to $d$. We have thus shown that the morphism from $V \oslash_Z S$ to $Z$ is generically of dimension $d$.
\end{dem}

\begin{prop}\label{PROP - point ab par enjoliveur}
Let $(X,x)$ be a strict pointed $K$-analytic space, $((Z,\varphi),S)$ a local enjoliveur of $(X,x)$ and $t$ a point in $\Abhy{Z}$, the point $\varphi(t)$ belongs to $\Abhy{X}$.
\end{prop}
\begin{dem}
We reduce to the case where $\varphi$ is either a morphism of $\fC$ or a blow-up.
\begin{enumerate}
\item If $\varphi$ is a morphism of $\fC$, it is, in particular, flat and of relative dimension zero. It follows that the dimension of $X$ at $\varphi(t)$ is equal to the dimension of $Z$ at $t$. Furthermore, the finiteness of $\CrpResC{t}$ over $\CrpResCR{\varphi(t)}$ ensures that:
\[ d_K(\varphi(t)) =  d_K(t) \]
\item If $\varphi$ is a blow-up, let us choose a centre $S$. The assumption made about $t$ implies that it does not lie on any Zariski closed set of $Z$ with empty interior. In particular, the point $t$ does not lie on the divisor $\BLupD{S}{X}$. The lemma then follows from the fact that $\varphi$ induces an isomorphism from $Z \priv \BLupD{S}{X}$ to $X \priv S$.
\end{enumerate}
\end{dem}

\subsection{Lifting local enjoliveurs}

\begin{prop}\label{PROP - remonter des enjoliveurs locaux}
Let $(X,x)$ be a strict pointed $K$-analytic space, $Y$ a closed $K$-analytic subspace of $X$ containing $x$, and $(W,T)$ a local enjoliveur of class $n$ of $(Y,x)$, then there exists a local enjoliveur $(Z,S)$ of class less than $n$ of $(X,x)$ such that: 
\begin{enumerate}
\item The $Y$-analytic space $W$ is equal to $(Y \times_X Z) \oslash_Z S$.
\item The divisor $T$ is equal to $W \times_Z S$.
\item The image of $S$ in $X$ is equal to that of $T$.
\end{enumerate}
\end{prop}
\begin{dem}
We prove this by recurrence with respect to $n$. If $n$ is zero, then $(W,T)$ is equal to $(Y,\emptyset)$ and we may take $(Z,S)$ to be the pair $(X,\emptyset)$. Now suppose that $n$ is non-zero and assume the result holds for integers strictly less than $n$. Since $(W,T)$ is of class $n$, there exists a local enjoliveur $(U,\Sigma)$ of $(Y,x)$ of class $n-1$, a finite family $(V_i)_{i \in I}$ of affinoid domains of $U$ and, denoting by $V$ the disjoint sum of $(V_i)_{i \in I} $, a closed analytic subspace $E$ of $V$ containing $V \times_U \Sigma$ such that $W$ is identified with the splitting of $V$ along $E$ and $T$ with its exceptional divisor.
\[ \xymatrix{ W = \BLup{E}{V} \ar[rr] &&  V = \disjsum_{i \in I} V_i \ar[rr] && U \ar[rr] && \ldots \ar[rr] && X } \]
The recurrence hypothesis then guarantees the existence of a local enjoliveur $(D,R)$ of $(X,x)$ of class less than $n-1$ such that:
\begin{enumerate}
\item The $Y-$space $U$ is equal to $(Y \times_X D) \oslash_D R$.
\item The quotient $\Sigma$ is equal to $U \times_D R$.
\item The image of $R$ on $X$ is equal to that of $\Sigma$.
\end{enumerate}
It follows from lemma \ref{LEM - Relevement de domaines rationnels} and the Gerritzen–Grauert theorem that there exists a finite family $(D_i')_{i \in I}$ of affinoid domains of $D$ such that, denoting by $D'$ the disjoint sum of $(D_i')_{i \in I} $, the $K$-analytic space $V$ is identified with $D' \times_D U$. Let $\cI$ be the coherent sheaf of ideals of $\cOf{D}$ defining $R$, $\cJ$ the coherent sheaf of ideals of $\cOf{D'}$ defining the image of $E$ in $D'$, and $F$ the closed analytic subspace of $D'$ defined by the sheaf of ideals $(\cI\cOf{D'}) \cap \cJ$, let $Z$ denote the blow-up of $D'$ along $F$ and $S$ the exceptional divisor of this blow-up. By construction, the pair $(Z,S)$ is a local enjoliveur of $(X,x)$ of class less than $n$. Now, since $W$ is identified with the blow-up of $V$ along $E$, since $\Sigma$ is equal to $U\times_D R$ and since $E$ contains $V \times_U \Sigma$, the closed analytic subspace $V \times_{D'} F$ of $V$ is identified with $E$. Consequently, the $K$-analytic space $W$ is the strict transform of $V$ relative to the blow-up $Z$ of $D'$ along $F$, so that $T$ is equal to $W\times_Z S$ and that:
\[ W = (Z \times_{D'} V) \oslash_Z S = (Z \times_D U) \oslash_Z S = (Z \times_D ((Y \times_X D) \oslash_D R)) \oslash_Z S = (Z \times_X Y) \oslash_Z S \]
where the last equality follows from the fact that $S$, by definition, contains $Z \times_D R$. The equality between $T$ and $W\times_Z S$ ensures that the image of $T$ in $X$ is contained in that of $S$. And the image of $S$ on $X$ is contained in that of $F$, which is itself contained in the union of the image of $E$ and the image of $R$, and therefore in the union of the image of $E$ and the image of $\Sigma$. Since $E$ contains $V \times_U \Sigma$, it follows that the image of $S$ on $X$ is contained in the image of $E$, and hence in that of $T$. Consequently, the image of $S$ on $X$ is ultimately equal to that of $T$.
\end{dem}

\newpage \section{Flattening of coherent modules}

\begin{att}
Throughout this section, $K$ denotes a non-trivially valued field.
\end{att}

\subsection{Flattening of coherent modules in the quasi-smooth case}

The property discussed in this subsection is a special case of the general flattening result established by Antoine Ducros. It is proven and used in a crucial way in \cite{Duc20}, as well as in the proof of our previously stated result. For the reader’s convenience, we reproduce here the direct proof, which is already included in Antoine Ducros’s proof.

\subsection*{Fitting ideals}

In the proof of the flattening theorem established by Antoine Ducros in \cite{Duc20}, the concept of a Fitting ideal is crucial. Let us briefly review what it entails.\par\medskip

Let $A$ be a commutative ring, $M$ a finite-type $A$-module, and $r$ a positive integer or zero. Let us consider a presentation of $M$:
\[ \xymatrix{ A^{(I)} \ar[rr]^u && A^n \ar[rr] && M \ar[rr] && 0 } \]
The $r$-th Fitting ideal of $M$ is the ideal of $A$ generated by the minors of size $n - r$ of the matrix of $u$ in the canonical bases if $r$ is less than $n$, and the ideal $A$ if $r$ is strictly greater than $n$. This terminology is justified because we show that this ideal does not depend on the chosen presentation of $M$.\par\medskip

Let $X$ be a scheme (or a $K$-analytic space) and $\cF$ a quasi-coherent sheaf of finite type (or a coherent sheaf) on $X$. By gluing from the affine (or affinoid) case, we can define the $r$-th Fitting ideal $\cI$ of $\cF$ on $X$. It is a quasi-coherent (resp. coherent) sheaf of ideals, whose locus of zeros is exactly the set of points in $X$ at which the pointwise rank of $\cF$ is strictly greater than $r$.

\subsection*{Coefficient ideals}

We now set out to show the existence of a "coefficient ideal" associated with a closed analytic subspace of the source of a quasi-smooth and compact morphism with geometrically irreducible fibers. This result will play a crucial role in the construction of blow-ups with a flattening effect. Indeed, the centers of these blow-ups will all be defined by well-chosen coefficient ideals.

\begin{theo}[Theorem (3.9) of \cite{Duc20}]\label{3.9 de Duc20}
Let $X$ and $Y$ be two $K$-analytic spaces, $Z$ a closed $K$-analytic subspace of $Y$ and $f$ a quasi-smooth, compact morphism from $X$ to $Y$ with geometrically irreducible fibers.
\begin{enumerate}
\item The category of analytic spaces $T$ defined on a complete extension of $K$ and endowed with a $K$-morphism to $X$ such that the closed immersion of $Z \times_X T$ into $Y \times_X T$ is an isomorphism has a terminal object $S$.
\item The space $S$ is $K$-analytic and the $K$-morphism from $S$ to $X$ is a closed immersion.
\item The $S$-algebra commutes with changes of basis.
\end{enumerate}
\end{theo}

\begin{deff}[Coefficient space]
Let $X$ and $Y$ be two $K$-analytic spaces, $Z$ a closed $K$-analytic subspace of $Y$ and $f$ a quasi-smooth, compact morphism from $X$ to $Y$ with geometrically irreducible fibers, we call the space of $f$-coefficients of $Z$, and denote by $\IdCoef{Z/f}$, the final object of the category of analytic spaces $T$ defined over a complete extension of $K$ and equipped with a $K$-morphism to $X$ such that the closed immersion of $Z \times_X T$ into $Y \times_X T$ is an isomorphism.
\end{deff}

\begin{voc}[Coefficient ideals]
Let $X$ and $Y$ be two $K$-analytic spaces, $\cI$ a sheaf of ideals of $\cOf{Y}$, $Z$ the closed $K$-analytic subspace of $Y$ defined by $\cI$ and $f$ a quasi-smooth, compact morphism from $X$ to $Y$ with geometrically irreducible fibers, the sheaf of ideals defining the space of $f$-coefficients of $Z$ is called the $f$-coefficient ideal of $\cI$.
\end{voc}

\begin{prop}[Theorem (3.9) of \cite{Duc20}]
Let $X$ and $Y$ be two $K$-analytic spaces, $\cI$ a sheaf of ideals of $\cOf{Y}$ and $\cJ$ the ideal of $f$-coefficients of $\cI$, if $\cJ$ is invertible, then the coherent subsheaf $(\cJ\cOf{Y})^{-1} \cI$ of $\cK_Y$ is contained in $\cOf{Y}$ and its coefficient ideal is $\cOf{X}$.
\end{prop}

\subsection*{Flattening of coherent modules in the quasi-smooth case}

\begin{lem}\label{LEM - Aplatissement fort 1}
Let $X$ and $Y$ be two compact strict $K$-analytic spaces, $f$ a quasi-smooth morphism from $Y$ to $X$ with geometrically irreducible fibers of dimension $n$, $V$ a Zariski open set of $Y$ such that $\cF_V$ is $G$-locally free and $S$ a closed $K$-analytic subspace of $X$ such that $Y \times_X S$ is the union of the fibers of $f$ that do not intersect $V$, there exists a closed $K$-analytic subspace $S'$ of $X$ with support $S$ such that the coherent sheaf:
\[ (\cF_{Y \times_X \BLup{S}{X}}) \oslash_X S' \]
or $\BLup{S}{X}$ is flat in dimensions greater than $n$.
\end{lem}
\begin{dem}
Let $T$ be the complement of $V$ in $Y$ and $\Omega$ the open complement of $S$ in $X$, for every point $x$ in $\Omega$, the fiber $V_x$ is nonempty. The module $\cF$ is free at every point of $V$ and, in particular, at every point of $V_x$. Since the latter is a Zariski open set of the normal and connected space $Y_x$, it is connected by property (1.11) of \cite{Duc18} and the rank of $\cF$ is therefore equal to the same integer $r(x)$ at every point in $V_x$. Let $r$ denote the function:
\[ \left| \begin{array}{c c c}
\Omega & \to & \eN \\
x & \mapsto & r(x)
\end{array} \right. \]
\begin{enumerate}
\item Let $x$ be a point in $\Omega$ and let $\Omega'$ be a compact analytic neighborhood of $x$ in $\Omega$, for every $\omega$ in $\Omega'$, the closed set $T_\omega$ has an empty interior in $Y_\omega$. Corollary (3.8) of \cite{Duc18} then ensures that there exists a compact analytic domain $T'$ of $Y \times_X \Omega'$ that avoids $T \times_X \Omega'$ and intersects all fibers of $Y$ above $\Omega'$. Since $T'$ is contained in $V$, the sheaf $\cF_{T'}$ is $G$-locally free. For every integer $i$, let $T_i'$ denote the closed open set of $T'$ on which the rank of $\cF$ is equal to $i$. Then $T_i'$ is empty except for a finite number of integers and the family $(T_i')_{i \in \eN}$ forms an open covering of $T'$. Furthermore, for any $x$ in $\Omega$, since the rank is constant on $V_x$, the images of the nonempty open sets of the family $(T_i')_{i \in \eN}$ on $\Omega'$ are pairwise disjoint and compact. This gives us a representation of $\Omega'$ as a finite disjoint union of closed open sets on which $r$ is constant. In particular, the function $r$ is constant in the neighborhood of $x$.
\item Let $m$ be an integer, $I$ the set of integers $i$ such that $T_i'$ is nonempty and $\Omega_m$ the set of points $x$ in $\Omega$ such that $r(x) = m$. Since $\Omega_m$ is a closed open set in $\Omega$, it is a Zariski open set in $X$ by (1.12) of \cite {Duc18}. Then, by (5.11) in \cite{Duc18}, there exists a $K -$analytic closed subspace $S'$ of $X$ with support $S$ such that the finite partition $(\Omega_i \times_X \BLup{S'}{X})_{i \in I}$ of $\Omega \times_X \BLup{S'}{X}$ is induced by a finite partition into Zariski open sets of $\BLup{S}{X}$. By replacing $X$ with $\BLup{S}{X}$ and proceeding open by open, we reduce the problem to the case where $r$ is constant on $\Omega$, with value still denoted by $r$, and where $S$ is a divisor. Under these assumptions, since $\Omega$ is the complement of $S$, $\Omega$ is analytically dense, and every embedded component of $X$ intersects $\Omega$.
\item Let us show that $V$ is analytically dense in $Y$. Let $\Sigma$ be an assassin component of $Y$, since $f$ is quasi-smooth, the set $\Sigma$ is an irreducible component of $\reci{f}\fens{\Theta}$ for some embedded component $\Theta$ of $X$ equipped with its reduced structure. Since the restriction of $f$ to $\reci{f}\fens{\Theta}$ is flat, the image $f\fens{\Sigma}$ is Zariski-dense in $\Theta$, and, as an embedded component of $X$, the closed set $\Theta$ intersects $\Omega$, and thus the image $f\fens{\Sigma}$ is not contained in $S$ by what has been shown above. In other words, the set $\reci{f}\fens{S} \cap \Sigma$ has an empty interior in $\Sigma$. Consequently, there exists a point $y$ in $\Sigma$, whose image on $X$ is denoted by $x$, that lies neither on $\reci{f}\fens{S}$ nor on any irreducible component of $\reci{f}\fens{\Theta}$ other than $\Sigma$. The relative dimension of the restriction of $f$ to $\Sigma$ at $y$ is then equal to the relative dimension of $f$ at $y$, that is, to $n$. Consequently, since the fiber $Y_x$ is irreducible of dimension $n$, the set $\Sigma \cap Y_x$ is equal to $Y_x$. Furthermore, since $x$ is not in $S$, the fiber $V_x$ is nonempty, and $\Sigma$ therefore intersects $V$. Consequently, $V$ is analytically dense in $Y$.
\item Let $\cI$ be the $r$-th Fitting ideal of $\cF$ and $\cJ$ the ideal of the $f$-coefficients of $\cI$ given by (3.9) of \cite{Duc20}, since $\cF_V$ is G-locally free of rank $r$, we have that $\cI_V$ is equal to $\cOf{V}$, and since $V$ intersects all fibers of $f$ above $\Omega$, it follows that $\cJ_\Omega$ is equal to $\cOf{\Omega}$. Let $\cJ'$ be the invertible ideal sheaf defining $S$. By replacing $X$ with its blow-up along $\cJ \cJ'$, we reduce the problem to the case where $\cJ$ is invertible. Let $\cH$ be the coherent subsheaf $(\cJ\cOf{Y})^{-1}\cI$ of $\cK_Y$; it is in fact contained in $\cOf{Y}$, and its coefficient ideal is itself equal to $\cOf{X}$ by (3.9) of \cite{Duc20}. Let $Y'$ be the complementary open set of the Zariski closed set of $Y$ defined by the ideal sheaf $\cH$. By definition, $\cH_{Y'}$ is equal to $\cOf{Y'}$, and thus $\cI_{Y'}$ is equal to $\cJ\cOf{Y'}$. In particular, the ideal $\cI_{Y'}$ is invertible. Let $x$ be a point in $X$. Since the ideal of coefficients of $\cH$ is equal to $\cOf{X}$, the closed analytic subspace of $Y_x$ defined by $\cH_{Y_x}$ is not isomorphic to $Y_x$, which implies that its support is different from $Y_x$ because the space $Y_x$ is quasi-smooth, and therefore regular, and thus reduced. Consequently, the open set $Y'$ intersects all fibers of $f$. Note also that since $\cJ_\Omega$ is equal to $\cOf{\Omega}$, we have the equality:
\[ \cJ\cOf{Y \times_X \Omega} = \cOf{Y\times_X \Omega} \]
so we have the following equality:
\[ \cH_{Y \times_X \Omega} = \cI_{Y \times_X \Omega} \]
This results in the following equalities:
\[ \cH_V = \cI_V = \cOf{V} \quad\et\quad \cI_{Y' \times_X \Omega} = \cH_{Y' \times_X \Omega} = \cOf{ Y' \times_X \Omega} \]
and that $V$ is contained in $Y'$. We will now show that $\cF \oslash_X S$ is $X$-flat in dimensions greater than $n$. Let $\cN$ be the coherent subbundle of $\cF$ consisting of sections annihilated by $\cI$, since we have:
\[ \cI_{Y' \times_X \Omega} = \cOf{Y' \times_X \Omega} \]
The sheaf $\cN_{Y' \times_X\Omega}$ is zero. Consequently, the support of any section of $\cN$ over an open set of $Y'$ is contained in $Y' \times_X S$, and we therefore have an exact sequence of coherent sheaves on $Y'$:
\[ \xymatrix{0 \ar[rr] && \cP \ar[rr] && (\cF/\cN)_{Y'} \ar[rr] && (\cF \oslash_X S)_{Y'} \ar[rr] && 0 } \]
where $\cP$ is the coherent sheaf on $Y'$ defined by this exact sequence. We will first show that $(\cF/\cN)_{Y'}$ is G-locally free of rank $r$, since $f$ is quasi-smooth, this will ensure that $(\cF/\cN)_{Y'}$ is $X$-flat. Let $W$ be an affinoid domain of $Y$. Since the Zariski open set $V$ is analytically dense in $Y$, the intersection $V \cap W$ is analytically dense in $W$ by (2.16) of \cite{Duc20}. In other words, every assassin component of $W$ intersects $V$, which means, by (2.9.(2)) of \cite{Duc20}, that the assassin of $\Algbr{W}$ is contained in $\Algbr{(W \cap V)}$. Since the sheaf $\cF_V$ is G-locally free of rank $r$, the restriction of the sheaf $\Algbr{(\cF_W)}$ to $\Algbr{(W \cap V)}$ is locally free of rank $r$. Consequently, the sheaf $\Algbr{(\cF_W)}$ is locally free of rank $r$ at every point of the killer set of $\Algbr{W}$. It then follows from (5.4.3) of \cite{RG71} that $\Algbr{(\cF_W)}/\Algbr{(\cN_W)}$ is locally free of rank $r$ over the largest open set of $\Algbr{W}$ on which the $r$-th Fitting ideal $\Algbr{\cI}$ of $\Algbr{\cF}$ is invertible, an open set that contains $\Algbr{(Y' \cap W)}$. Since this holds for every affinoid domain $W$ of $Y$, the sheaf $(\cF/\cN)_{Y'}$ is G-locally free of rank $r$ and therefore $X$-flat. Since $V$ does not intersect $Y \times_X S$, we have the following equalities:
\[ (\cF/\cN)_V = (\cF \oslash_X S)_V = \cF_V \]
and the sheaf $\cP_V$ is zero. Furthermore, since the coherent sheaf $\cP$ is a subsheaf of the G-locally free sheaf $(\cF/\cN)_{Y'}$, any embedded component of $\cP$ is then a embedded component of $Y'$, and therefore intersects $V$. But since $\cP_V$ is zero, this implies that $\cP$ has no embedded component, which means it is zero. Consequently, the sheaf $(\cF \oslash_X S)_{Y'}$ is equal to $(\cF/\cN)_{Y'}$, and $(\cF \oslash_X S)_{Y'}$ is therefore $X$-flat. Since $Y'$ intersects all fibers of $f$, the coherent sheaf $\cF \oslash_X S$ is $X$-flat in dimensions greater than $n$.
\end{enumerate}
\end{dem}

\begin{prop}\label{PROP - Aplatissement fort sur quasi-lisse}
Let $X$ and $Y$ be two compact strict $K$-analytic spaces and let $f$ be a quasi-smooth morphism from $Y$ to $X$ with geometrically irreducible fibers of dimension $n$, there exists a closed $K$-analytic subspace $T$ of $Y$ containing $\LNPlat[\geqslant n]{\cF}{X}$ such that:
\begin{enumerate}
\item The closed set $T$ has an $f$-coefficient space $S$ whose support is $f\fens{\LNPlat[\geqslant n]{\cF}{X}}$ and whose coherent sheaf is:
\[ (\cF_{Y \times_X \BLup{S}{X}}) \oslash_X S \]
is $\BLup{S}{X}-$flat in dimensions greater than $n$.
\item For every point $s$ in $S$, there exists a point $t$ in $T_s \cap \LNPlat[\geqslant n]{\cF}{X}$ such that:
\[ d_{\CrpResC{s}}(t) = n \]
\end{enumerate}
\end{prop}
\begin{dem}
Under our assumptions, the set $\LNPlat[\geqslant n]{\cF}{X}$ is exactly the union of the fibers of $f$ that do not intersect $\LPlat{\cF}{X}$. Let $U$ be the complement of $\LNPlat[\geqslant n]{\cF}{X}$ and $V$ be the set of points where $\cF$ is free, these are two Zariski open sets of $Y$. Since $Y$ is quasi-smooth over $X$, the coherent sheaf $\cOf{Y}$ is flat over $X$, and $V$ is therefore contained in $\LPlat{\cF}{X}$, and a fortiori in $U$.  Let $x$ be a point in $X$ such that $U_x$ is nonempty, this means that $Y_x$ intersects $\LPlat{\cF}{X}$, and we then deduce from (4.13) in \cite{Duc20} that $V$ intersects $Y_x$. Since this holds for every $x$ such that $U_x$ is nonempty, we can also describe $\LNPlat[\geqslant n]{\cF}{X}$ as the union of the fibers of $f$ that do not intersect $V$. Let $T$ be the complement of $V$ in $Y$ equipped with its reduced structure, where the fibers of $f$ are geometrically irreducible, the closed set $T$ in $Y$ has a coefficient ideal on $X$ by (3.9) of \cite{Duc20}, defining a closed $K$-analytic subspace $S$ of $X$. By definition, a point $x$ in $X$ belongs to $S$ if and only if its fiber $Y_x$ does not intersect $V$, that is, if and only if $x$ belongs to the image of $\LNPlat[\geqslant n]{\cF}{X}$ under $f$: indeed, since the fiber $Y_x$ is quasi-smooth, it is regular and, in particular, reduced. Thus, a closed $K$-analytic subspace of $Y_x$ is isomorphic to $Y_x$ if and only if its support is equal to $Y_x$. In other words, the support of $S$ is equal to the image of $\LNPlat[\geqslant n]{\cF}{X}$ under $f$, and the result then follows from the lemma \ref{LEM - Aplatissement fort 1}.
\end{dem}

\subsection{Flattening of coherent modules in the boundaryless case}

\begin{lem}\label{LEM - Aplatissement fort 2}
Let $(X,x)$ and $(Y,y)$ be two proper, strict pointed $K$-analytic spaces, $f$ be a pointed morphism from $(Y,y)$ to $(X,x)$ such that $y$ lies in the interior of $f$, $S$ a divisor of $X$, $n$ a non-zero integer and $\cF$ a coherent $\cOf{Y}$-module such that $\cF \oslash_X S$ is flat in dimensions greater than $n+1$, there exists a strict affinoid neighbourhood $W$ of $x$ in $X$, a closed $K$-analytic subspace $\Sigma$ of $W$ and an analytic neighbourhood $U$ of $y$ in $Y$ such that, denoting by $V$ the $K$-analytic space $U \times_X \BLup{\Sigma}{W}$, we have:
\begin{enumerate}
\item The module $\cF_V \oslash_W \Sigma$ is $\BLup{\Sigma}{W}-$flat in dimensions greater than $n$.
\item The $K$-analytic space $\Sigma$ is set-theoretically equal to:
\[ f\fens{\LNPlat[\geqslant n]{\cF_{U \times_X W}}{W}} \cup (S \cap W) \]
\item The closed set $V \times_X \Sigma$ of $V$ is contained in the inverse image of:
\[ (U \times_X (S \cap W)) \cup \LNPlatS[\geqslant n]{\cF_{U \times_X W}}{W} \]
\item The image of $V$ in $U \times_X W$ contains the complement of:
\[ (U \times_X (S \cap W)) \cup \LNPlatS[\geqslant n]{\cF_{U \times_X W}}{W} \]
\item The image of $V \times_X \Sigma$ in $U \times_X W$ contains the set of points in:
\[ (U \times_X (S \cap W)) \cup \LNPlat[\geqslant n]{\cF_{U \times_X W}}{W} \]
whose image on $W$ is not adherent to:
\[ \left( f\fens{\LNPlat[\geqslant n]{\cF_{U \times_X W}}{W}} \cup ( f\fens{U \times_X W} \cap S ) \right) \cap \Abhy{X} \]
\item For any point $s$ in $\Sigma \priv (S \cap W)$, there exists a point $u$ in:
\[ U_s \cap \LNPlat[\geqslant n]{\cF_{U \times_X W}}{W} \]
such that $d_{\CrpResC{s}}(u) = n$.
\end{enumerate}
\end{lem}
\begin{dem}
To establish the desired property, we may reason locally on $X$ around $x$ and on $Y$ around $y$. We may therefore assume that both $X$ and $Y$ are strict affinoid. Furthermore, in this demonstration, we shall allow ourselves to restrict $X$ around $x$ and $Y$ around $y$ as much as necessary. By reasoning locally one last time on $Y$ around $y$, we may, by the property \ref{THEO - existence de devissage connexe}, assume that we are in one of the following two cases:
\begin{itemize}
\item \textit{The point $y$ does not belong\footnote{In particular, if the point $y$ does not belong to the support of $\cF \oslash_X S$.} to $\LNPlat[\geqslant n]{\cF}{X}$.}\\
We can then restrict $Y$ around $y$ to reduce the problem to the case where $\cF \oslash_X S$ is $X$-flat in dimensions greater than $n$, and the desired assertion holds for $W = X$, with $\Sigma = S$ and\footnote{Now, since $S$ is a divisor, the blow-up $\BLup{S}{X}$ is isomorphic to $X$ and the result is the empty.} $U = Y$.
\item \textit{The point $y$ belong to $\LNPlat[\geqslant n]{\cF}{X}$.} \\
Even if we restrict $Y$ to the neighbourhood of $y$, there exists a total Raynaud–Gruson's $X$-dévissage:
\[ (Y,\{T_i,t_i,\pi_i,u_i,\cL_i,\cP_i\}_{1 \leqslant i \leqslant r}) \]
of $\cF \oslash_X S$ at $y$ in dimension $(n_i)_{1 \leqslant i \leqslant r}$ bounded by the dimension of $\SuppFa{\cF \oslash_X S}_x$. Even if this means shortening the decomposition, which may then no longer be total, we can replace the condition $\cP_r = 0$ with the following two conditions:
\begin{enumerate}[label=\arabic*)]
\item The dimension of $\SuppFa{\cP_r}_x$ is strictly less than $n$. \label{LEM - Aplatissement 2 cdt 1}
\item For every integer $i$ in $\lseg 1,r-1 \rseg$, the dimension of $\SuppFa{\cP_i}_x$ is greater than $n$. \label{LEM - Aplatissement 2 cdt 2}
\end{enumerate}
For every integer $i$ in $\lseg 1,r-1 \rseg$, this therefore implies that $n_i$ is greater than $n+1$ and that $n_{i+1}$ is greater than $n$. By (4.12) in \cite{Duc20}, we may also assume, even if this requires restricting $X$, that the following two statements hold:
\begin{enumerate}[label=\arabic*),start=3]
\item For every integer $i$ in $\lseg 1,r \rseg$ and every $\xi$ in $X$, the intersection $\SuppFa{\cP_i} \cap (T_i)_\xi$ has an empty interior in $(T_i)_\xi$.
\item For every integer $i$ in $\lseg 1,r \rseg$ such that $n_i$ is greater than $n+1$, the morphism $u_i$ is universally $X$-injective.
\end{enumerate}
We are then in one of the following three situations:
\begin{enumerate}
\item \textit{Let us assume that $n_r$ is strictly less than $n$.}\\
Since, for every $i$ between $1$ and $r-1$, the integer $n_{i+1}$ is greater than $n$ by condition \ref{LEM - Aplatissement 2 cdt 2}, it follows that $r$ is equal to $1$. The relative dimension of the support of $\cF \oslash_X S$ on $X$ is therefore strictly less than $n$, and $\cF \oslash_X S$ is, in particular, $X$-flat in dimensions greater than $n$. The desired assertion is therefore true with $W = X$ and $\Sigma = S$.
\item \textit{Let us assume that $n_r$ is strictly greater than $n$.}\\
In this case, for every $i$ between $1$ and $r$, the morphism $u_i$ is universally injective and, together with the condition \ref {LEM - Aplatissement 2 cdt 1}, this implies, by (4.12) of \cite{Duc20}, that $\cF \oslash_X S$ is $X$-flat in dimensions greater than $n$, and the desired assertion therefore holds with $W = X$ and $\Sigma = S$.
\item \textit{Let us assume that $n_r$ is equal to $n$.}\\
Let $W$ denote the image of $T_r$ in $X$ and $\omega$ the morphism from $T_r$ to $W$. Since $T_r$ is compact and $t_r$ lies in the interior of $\omega$, the morphism $\omega$ is open on a neighbourhood of $t_r$, and $W$ is a neighbourhood of $x$ in $X$. Even if we restrict $X$ to the neighbourhood of $x$, we can therefore assume that $\omega$ is surjective.
\begin{itemize}
\item Since $\omega$ is quasi-smooth with geometrically connected fibres and surjective, the morphism from $T_r$ to $X$ is quasi-smooth with geometrically irreducible fibres of dimension $n$. Let $\cQ$ denote the sheaf $(\pi_r)_\ast(\cP_{r-1})$, by the property \ref{PROP - Aplatissement fort sur quasi-lisse}, there exists a closed analytic subspace $\Theta$ of $T_r$ with an $\omega$-coefficient space $F$ whose support is the image of $\LNPlat[\geqslant n]{\cQ}{X} = \LNPlatS[\geqslant n]{\cQ}{X}$ in $X$, and such that the sheaf:
\[ \cQ_{T_r \times_X \BLup{F}{X} } \oslash_X F \]
Let $\BLup{F}{X}$ be flat in dimensions greater than $n$. Let $\Sigma$ then be the closed $K$-analytic subspace of $X$ defined by the product of the ideals defining $F$ and $S$ respectively, its support is the union of $S$ and the image of $\LNPlat[\geqslant n]{\cQ}{X}$, and the blow-up $\BLup{\Sigma}{X}$ is equal to $\BLup{F}{X}$ by (5.5.6) in \cite{Duc20}. Furthermore, the locus of $\BLup{\Sigma}{X}$-flatness of the sheaf:
\[ \cQ_{T_r \times_X \BLup{F}{X} } \oslash_X \Sigma \]
contains that of the sheaf:
\[ \cQ_{T_r \times_X \BLup{\Sigma}{X} } \oslash_X F \]
Indeed, the first is the quotient of the second by its subsheaf consisting of sections with support in $T_r \times_X \BLupD{\Sigma}{X}$, which is a divisor of $T_r \times_X \BLup{\Sigma}{X}$ since the morphism from $T_r$ to $X$ is flat. Our statement then follows from point (6.4) of \cite{Duc20}. Consequently, the sheaf:
\[ \cQ_{T_r \times_X \BLup{F}{X} } \oslash_X \Sigma \]
is $\BLup{\Sigma}{X}-$flat in dimensions greater than $n$.
\end{itemize}
We shall now prove two statements:
\begin{enumerate}[label=(\alph*)]
\item the coherent sheaf : \label{LEM - Aplatissement fort 2 cdtB}
\[ (\cF \oslash_X S)_{Y \times_X \BLup{\Sigma}{X}} \oslash_X \Sigma = \cF_{Y \times_X \BLup{\Sigma}{X}} \oslash_X \Sigma \]
is $\BLup{\Sigma}{X}-$flat in dimension greater than $n$.\\
Let $\cG$ denote the sheaf $(\cF \oslash_X S)_Y$. By changing the basis via the morphism from $\BLup{F}{X}$ to $X$, we get, without altering the properties of $\cD$, at the case where $\Sigma$ is a divisor, where $\BLup{\Sigma}{X}$ is equal to $X$ and where $\cQ \oslash_X \Sigma$ is $X$-flat in dimensions greater than $n$. It suffices to show that this implies the $X$-flatness of $\cG \oslash_X \Sigma$ in dimensions greater than $n$. Since $\cQ \oslash_X \Sigma$ is $X$-flat in dimensions greater than $n$, the same holds for $\cP_{r-1} \oslash_X \Sigma$ by (4.1.15 (3)) of \cite{Duc18}. For every $i$ between $1$ and $r-1$, the exact sequence:
\[ \xymatrix{ 0 \ar[r] & \cL_i \ar[rr] && (\pi_i)_\ast(\cP_{i-1}) \ar[rr] && \cP_i \ar[r] & 0 } \]
By setting $T_0 = Y$ and $\cP_0 = \cG$, this implies, by Lemma 5.5.3 (ii) of \cite{RG71}, an exact sequence:
\[ \xymatrix{ 0 \ar[r] & \cL_i \ar[rr] && (\pi_i)_\ast(\cP_{i-1}) \oslash_X \Sigma \ar[rr] && \cP_i \oslash_X \Sigma \ar[r] & 0 } \]
hence a Raynaud-Gruson's dévissage $\cD'$ of $\cG \oslash_X \Sigma$ in dimensions $(n_i)_{1 \leqslant i \leqslant r-1}$. Since the linear maps constituting it are injective, and since the coherent sheaf $\cP_{r-1} \oslash_X \Sigma$ is $X -$flat in dimensions greater than $n$, it follows from Proposition (4.12) of \cite{Duc20} that the coherent sheaf $\cG \oslash_X \Sigma$ is $X-$flat in dimensions greater than $n$.
\item We have the equality : \label{LEM - Aplatissement fort 2 cdtA}
\[ f\fens{\LNPlat[\geqslant n]{\cF_Y}{X}} \cup S = \omega\fens{\LNPlat[\geqslant n]{\cQ}{X}} \cup S \]
Let $\xi$ be a point in $X$, the dimension of the support of $(\cP_r)_\xi$ is strictly less than $n_r = n$, and we know, moreover, that $u_i$ is universally $X$-injective for $i$ ranging from $1$ to $r-1$. By applying Proposition (4.12) of \cite{Duc20} on the one hand to the Raynaud–Gruson's $X$-dévissage of a single step $\{ T_r, \Id{}, u_r, \cL_r, \cP_r \}$ of $\cQ$ and, on the other hand, to $\cD$, we see that the following statements are equivalent:
\begin{itemize}
\item The dimension of $\LNPlat[\xi]{\cQ}{X}$ is strictly less than $n$.
\item The morphism $u_r$ is injective at every point of $(T_r)_\xi$.
\item The dimension of $\LNPlat[\xi]{(\cF \oslash_X S)_Y}{X}$ is strictly less than $n$.
\end{itemize}
If $\xi$ does not belong to $S$, the third statement simply amounts to requiring that the dimension of $\LNPlat{\cF_Y}{X}$ be strictly less than $n$. The equivalence between the first and third statements therefore implies that:
\[ f\fens{\LNPlat[\geqslant n]{\cF_Y}{X}} \priv S = \omega\fens{\LNPlat[\geqslant n]{\cQ}{X}} \priv S \]
hence the assertion.
\end{enumerate}
These two statements imply that $f$ satisfies the desired property, provided that the analytic domain of $Y \times_X \BLup{\Sigma}{X}$ is taken to be the entire space.
\begin{itemize}
\item The first assertion is a direct consequence of\ref{LEM - Aplatissement fort 2 cdtB}.
\item The fourth and fifth assertions follow from \ref{LEM - Aplatissement fort 2 cdtA}, from the fact that $\Sigma$ is set-theoretically equal to:
\[ \omega\fens{\LNPlat[\geqslant n]{\cQ}{X}} \cup S \]
and because the blow-up morphism from $\BLup{\Sigma}{X}$ to $X$ induces an isomorphism over $X \priv \Sigma$, which implies that the morphism from $Y \times_X \BLup{\Sigma}{X}$ to $Y \times_X X$ induces an isomorphism over $Y \priv (Y \times_X \Sigma)$.
\item Let $y$ be a point on $(Y \times_X S) \cup \LNPlatS[\geqslant n]{\cF_Y}{X}$ such that its image $x$ under $f$ is not adherent to:
\[ (f\fens{\LNPlatS[\geqslant n]{\cF_Y}{X}} \cup (f\fens{Y} \cap S)) \cap \Abhy{X} = (\omega\fens{\LNPlat[\geqslant n]{\cQ}{X}} \cup (f\fens{Y} \cap S)) \cap \Abhy{X} \]
Let us show that $y$ lies in the image of $Y \times_X \BLup{\Sigma}{X}$. The point $x$ lies in $\Sigma$ by definition, and no irreducible component of $X$ contained in $\Sigma$ passes through $x$. Indeed, if this were the case, that component would be contained in:
\[ f\fens{\LNPlatS[\geqslant n]{\cF_Y}{X}} \]
since $S$ is a divisor of $X$, and in particular has an empty interior in $X$, and $x$ would be a member of:
\[ f\fens{\LNPlatS[\geqslant n]{\cF_Y}{X}} \cap \Abhy{X} \]
which contradicts our hypothesis. Since no irreducible component of $X$ contained in $\Sigma$ passes through $x$, the point $x$ belongs to the image of the morphism from $\BLup{\Sigma}{X}$ to $X$ by (5.9.2) of \cite{Duc20}. It follows that $y$ belongs to the image of $Y \times_X \BLup{\Sigma}{X}$.
\item Let $s$ be a point in $\Sigma$ not belonging to $S$. Then the point $s$ belongs to $\omega\fens{\LNPlat[\geqslant n]{\cQ}{X}}$, and there exists a point $t$ in $\Theta_s \cap \LNPlat[\geqslant n]{\cQ}{X}$ such that:
\[ d_{\CrpResC{s}}(t) = n \]
Since we have :
\[ f\fens{\LNPlat[\geqslant n]{\cF_Y}{X}} \priv S = \omega\fens{\LNPlat[\geqslant n]{\cQ}{X}} \priv S \]
The point $s$ is the image of a point $u$ belonging to $\LNPlat[\geqslant n]{\cF_Y}{X}$ and $t$ is the image of $u$ under a composition of finite morphisms, so that we have:
\[ d_K(t) = d_K(u) = n \]
\end{itemize}
\end{enumerate}
\end{itemize}
\end{dem}

\begin{theo}[Flattening of coherent modules with local enjoliveur]\label{THEO - Applatissement fort}
Let $(X,x)$ and $(Y,y)$ be two good, strict pointed $K$-analytic spaces, $f$ be a pointed morphism from $(Y,y)$ to $(X,x)$ such that $y$ lies in the interior of $f$, $N$ and $n$ be two integers such that $n$ is less than $N$ and $\cF$ be a $\cOf{Y}$-module that is flat and coherent in dimensions greater than $N$, there exists a local enjoliveur $(Z,S)$ of $(X,x)$ of class less than $N - n$ and an analytic neighbourhood $U$ of $y$ in $Y$ such that, denoting by $V$ the $K$-analytic space $U \times_X Z$, we have:
\begin{enumerate}
\item The module $\cF_V \oslash_Z S$ is $Z$-flat in dimensions greater than $n$.
\item The set $V \priv (V \times_Z S)$ is contained in the inverse image of $Y \priv \LNPlat[\geqslant n]{\cF}{X}$.
\item The closed set $S$ of $Z$ is contained in the inverse image of $f\fens{\LNPlat[\geqslant n]{\cF_U}{X}}$.
\item The set $V \times_Z S$ is contained in the inverse image of $\LNPlatS[\geqslant n]{\cF}{X}$.
\item The image of $V$ in $U$ contains $U \priv \LNPlatS[\geqslant n]{\cF_U}{X}$.
\item The image of $V \times_Z S$ in $U$ contains every point in $\LNPlat[\geqslant n]{\cF_U}{X}$ whose image on $X$ is not adherent to:
\[ f\fens{ \LNPlat[\geqslant n]{\cF_U}{X} } \cap \Abhy{X} \]
\item The image of $S$ on $X$ is contained in $f\fens{\SuppFa{\cF_U}}$ and, if $\cF_U$ is non-zero, this image has a dimension less than:
\[ \dimAn{K}{\SuppFa{\cF_U}} - n \]
\end{enumerate}
\end{theo}
\begin{dem}\sautlistebis
\begin{enumerate}[label=(\Alph*)]
\item \textit{Automatic nature of the second assertion:}\label{THEO - Aplatissement fort cdt auto} \\
Let assume the first assertion. Let $t$ be a point in $V$ that does not lie above $S$ and $y_t$ be the image of $t$ on $Y$, the first assertion ensures that $t$ does not belong to $\LNPlat[\geqslant n]{\cF}{X}$ and, since $t$ is not situated above $S$, its image on $Z$ belongs to $\LPlat{Z}{X}$. It then follows by the flat descent of flatness that $y$ does not belong to $\LNPlat[\geqslant n]{\cF}{X}$.

\item \textit{Proof of the six other statements of the theorem:}\label{THEO - Aplatissement fort 1-6} \\
We shall prove the general case by descending induction on the integer $n$ less than $N$ for the first fifth other statements of the theorem, and by strong descending induction on the integer $n$ less than $N$ for the seventh statement. We therefore assume that the theorem holds for every integer strictly greater than $n$ and less than or equal to $N$.
\begin{enumerate}
\item The case where $n$ equals $N$: \\
Since $\cF$ is $X$-flat in dimensions higher than $N$, the conclusions of the theorem hold for the local enjoliveur $(X,\emptyset)$ of $(X,x)$ and the neighbourhood $Y$ of $y$.
\item The case where $n$ is strictly less than $N$: \\
The recurrence hypothesis guarantees that there exists a local enjoliveur $((Z,\varphi),S)$ of $(X,x)$ of class less than $N - (n+1)$ given by a factorization:
\[ \varphi : \xymatrix{ Z = Z_{2(N-(n+1))} \ar[r] & Z_{2(N-(n+1))-1} \ar[r] & \ldots \ar[r] & Z_0 = X } \]
and, for every $i$ between $0$ and $2(N-(n+1))$, a closed $K$-analytic subspace $S_i$ of $Z_i$ such that:
\[ S_{2(N-(n+1))} = S \quad\et\quad S_0 = \emptyset \]
Furthermore, the recurrence hypothesis guarantees that there exists an analytic neighborhood $U$ of $y$ in $Y$ such that, letting $V$ denote the $K$-analytic space $U \times_X Z$, we have:
\begin{enumerate}[label=(\arabic*${}^\ast$)]
\item The coherent sheaf $\cF_V \oslash_Z S$ is $Z$-flat in dimensions greater than $n+1$.
\item The open set $V \priv (V \times_Z S)$ of $V$ is contained in the inverse image of $Y \priv \LNPlatS[\geqslant n+1]{\cF}{X}$. \label{THEO - Aplatissement fort cdt2*}
\item The closed set $S$ of $Z$ is contained in the inverse image of $f\fens{\LNPlat[\geqslant n+1]{\cF_U}{X}}$. \label{THEO - Aplatissement fort cdt3*}
\item The closed set $V \times_Z S$ in $V$ is contained in the inverse image of $\LNPlatS[\geqslant n+1]{\cF}{X}$. \label{THEO - Aplatissement fort cdt4*}
\item The image of $V$ in $U$ contains $U \priv \LNPlatS[\geqslant n+1]{\cF_U}{X}$. \label{THEO - Aplatissement fort cdt5*}
\item The image of $V \times_Z S$ in $U$ contains every point in $\LNPlat[\geqslant n+1]{\cF_U}{X}$ whose image on $X$ is not adherent to:
\[ f\fens{\LNPlat[\geqslant n]{\cF_U}{X}} \cap \Abhy{X} \] \label{THEO - Aplatissement fort cdt6*}
\item If $N-(n+2)$ is positive, then for every $i$ between $0$ and $N-(n+2)$, the image of $S_{2i}$ on $X$ is contained in $f\fens{\SuppFa{\cF}}$, and if $\cF$ is nonzero, this image has a dimension less than:
\[ \dimAn{K}{\SuppFa{\cF}} - (N - (i - 1)) \]\label{THEO - Aplatissement fort cdt7*}
\end{enumerate}
For every point $t$ in $V$ whose image in $Y$ is $y$ and whose image in $Z$ is denoted by $z_t$, there exists, by lemma \ref{LEM - Aplatissement fort 2}, a strictly affinoid neighborhood $W_t$ of $z_t$ in $Z$, a closed $K$-analytic subspace $\Sigma_t$ of $W_t$ containing $S \cap W_t$, and an analytic neighborhood $U_t$ of $t$ in $V$ such that, denoting by $V_t$ the $K -$analytic space $U_t \times_Z \BLup{\Sigma_t}{W_t}$ and $\omega$ the morphism from $V$ to $Z$, we have the following commutative diagram:
\[ \xymatrix{
				&& (U_t,t) \ar[d] 							&& && V_t \cartesien  \ar[llll] \ar[ddd] \\
(U,y) \ar[d] 	&& (V,t) \cartesien \ar[dd]^\omega \ar[ll]_{\varphi_V} 	&& \\
(Y,y) \ar[d]_f 	&& 											&& \\
(X,x) 			&& (Z,z_t) \ar[ll]_\varphi && (W_t,z_t) \ar[ll] && \BLup{\Sigma_t}{W_t} \ar[ll] \\
} \]
and the following properties:
\begin{enumerate}[label=(\alph*)]
\item The module :
\[ \cF_{V_t} \oslash_{W_t} \Sigma = (\cF_V \oslash_Z S)_{V_t} \oslash_{V_t} \Sigma  \]
is $\BLup{\Sigma_t}{W_t}-$flat in dimension greater than $n$.\label{THEO - Aplatissement fort cdtA}
\skipitems{1}
\item The $K$-analytic space $\Sigma$ is, set-theoretically, equal to:
\[ \omega\fens{\LNPlat[\geqslant n]{\cF_{U_t \times_Z W_t}}{W_t}} \cup (S \cap W_t)  \]
\item The closed set $V_t \times_Z \Sigma$ of $V_t$ is contained in the inverse image of:
\[ (U_t \times_Z (S \cap W_t)) \cup \LNPlatS[\geqslant n]{\cF_{U_t \times_Z W_t}}{W_t} \]\label{THEO - Aplatissement fort cdtD}
\item The image of $V_t$ in $V \times_Z W_t$ contains the complement of:
\[ (U_t \times_Z (S \cap W_t)) \cup \LNPlatS[\geqslant n]{\cF_{U_t \times_Z W_t}}{W_t} \]\label{THEO - Aplatissement fort cdtE}
\item The image of $V_t \times_Z \Sigma$ in $U_t \times_Z W_t$ contains the set of points in:
\[ (U_t \times_Z (S \cap W_t)) \cup \LNPlatS[\geqslant n]{\cF_{U_t \times_Z W_t}}{W_t} \]
whose image on $W_t$ is not in the closure of :
\[ \left( \omega\fens{\LNPlat[\geqslant n]{\cF_{U_t \times_Z W_t}}{W_t}} \cup \left( \omega\fens{U_t \times_Z W_t} \cap S \right) \right) \cap \Abhy{Z} \]\label{THEO - Aplatissement fort cdtF}
\item For every point $s$ in $\Sigma \priv (S \cap W_t)$, there exists a point $u$ in:
\[ (U_t \times_Z W_t)_s \cap \LNPlat[\geqslant n]{\cF_{U_t \times_Z W_t}}{W_t} \]\label{THEO - Aplatissement fort cdtG}
such that $d_{\CrpResC{s}}(u) = n$.
\end{enumerate}
Since the morphism $\varphi$ is a composition of domanial morphisms and blow-ups, it is a topologically proper morphism. Hence, by change of basis, $\varphi_V$ is topologically proper. The fiber $V_y$ of $\varphi_V$ is then compact, and thus there exists a finite set of points $\Delta$ in $V_y$ such that the interiors of the sets in the family $(W_t)_{t \in \Delta}$ form a covering of $V_y$. Let $T$ denote the union of $(V_t)_{t \in \Delta}$, the set $T$ is a neighborhood of $V_y$ in $V$, and since $\varphi_V$ is topologically proper, the image of $T$ under $\varphi_V$ is a neighborhood of $y$ in $Y$. Furthermore, without loss of generality, we may restrict $U$ to the neighborhood of $y$ and we can assume that $V$ is equal to $T$. Let $W$ denote the disjoint sum of $(W_t)_{t \in \Delta}$, $V_1$ the disjoint sum of $(V_t)_{t \in \Delta}$, and $\Sigma$ the disjoint sum of $(\Sigma_t)_{t \in \Delta}$. The pair:
\[  (\BLupD{\Sigma}{W} , \xymatrix{ \BLup{\Sigma}{W} \ar[rr] && X } ) \]
defines a local enjoliveur $(\BLup{\Sigma}{W},\BLupD{\Sigma}{W})$ of $(X,x)$ of class less than $N - n$, whose data, together with that of $U$, satisfies the conclusions of the theorem.

\begin{enumerate}[label=\arabic*.]
\item The first assertion follows from \ref{THEO - Aplatissement fort cdtA}.
\skipitems{1}

\item \textit{Verification of the third assertion:} \\
It follows from \ref{THEO - Aplatissement fort cdt3*}.

\item \textit{Verification of the fourth assertion:} \\
Let $u$ be a point in $U \priv \LNPlatS[\geqslant n]{\cF_U}{X}$. It belongs, a fortiori, to $U \priv \LNPlatS[\geqslant n+1]{\cF_U}{X}$ and is therefore equal, according to \ref{THEO - Aplatissement fort cdt5*}, to the image of a point $v$ in $V$. Since $u$ does not belong to $\LNPlatS[\geqslant n]{\cF}{X}$, condition \ref{THEO - Aplatissement fort cdt4*} implies that $v$ does not belong to $V \times_Z S$. It is therefore located above the locus of flatness of the morphism from $Z$ to $X$. Consequently, since $u$ does not belong to $\LNPlatS[\geqslant n]{\cF_U}{X}$, the point $v$ does not belong to:
\[ \LNPlatS[\geqslant n]{\cF_V}{Z} \]
and therefore not to:
\[ \LNPlatS[\geqslant n]{\cF_V \oslash_Z S}{Z} \]
because $\cF_V \oslash_Z S$ coincides with $\cF_V$ above $Z \priv S$. It follows that, by \ref{THEO - Aplatissement fort cdtE}, $v$ belongs to the image of $V_1$ in $V$, and $u$ therefore belongs to the image of $V_1$ in $Y$.

\item \textit{Verification of the fifth assertion:} \\
Let $u$ be a point in $V_1 \times_W \Sigma$ and $\eta$ its image in $V$, it follows from \ref{THEO - Aplatissement fort cdtD} that $\eta$ belongs to:
\[ (V \times_Z (S \cap W)) \cup \LNPlatS[\geqslant n]{\cF_{V \times_Z W}}{W} \]
There are two cases to consider:
\begin{itemize}
\item If $\eta$ belongs to $(V \times_Z (S \cap W))$, then, according to \ref{THEO - Aplatissement fort cdt4*}, it lies above $\LNPlatS[\geqslant n+1]{\cF}{X}$, and a fortiori above $\LNPlatS[\geqslant n]{\cF}{X}$.
\item Otherwise, $\eta$ belongs to:
\[ \LNPlatS[\geqslant n]{\cF_{V \times_Z W}}{W} \]
and the same line of reasoning used in the proof of the fourth claim shows that $\eta$ lies above $\LNPlatS[\geqslant n]{\cF}{X}$.
\end{itemize}

\item \textit{Verification of the sixth assertion:} \\
Let $s$ be a point in $\LNPlat[\geqslant n]{\cF_U}{X}$ whose image $t$ on $X$ is not adherent to:
\[ f\fens{\LNPlat[\geqslant n]{\cF_U}{X}} \cap \Abhy{X} \]
and let $\Lambda$ be the set :
\[ (V \times_Z S) \cup \LNPlat[\geqslant n]{\cF_V \oslash_Z S}{Z} \]
Let's prove that $s$ is in the image of $V_1 \times_W \Sigma$, and to do so, let's first show that $s$ has a preimage in $\Lambda$.
\begin{itemize}
\item Suppose that $t$ does not belong to $f\fens{\LNPlat[\geqslant n+1]{\cF_U}{X}}$; then the point $s$ does not belong to $\LNPlat[\geqslant n+1]{\cF_U}{X}$. It follows, by \ref{THEO - Aplatissement fort cdt5*}, that there exists a point $v$ in $V$ lying above $s$. By \ref{THEO - Aplatissement fort cdt4*}, the point $v$ does not belong to $V \times_Z S$ and, consequently, it lies above the flattening locus of the morphism from $Z$ to $X$, so that it belongs to $\LNPlat[\geqslant n]{\cF_V}{Z}$, and therefore to:
\[ \LNPlat[\geqslant n]{\cF_V \oslash_Z S}{Z} \]
since it is not located above $S$. It therefore does indeed belong to $\Lambda$.
\item Let us now assume that $t$ belongs to $f\fens{\LNPlat[\geqslant n+1]{\cF_U}{X}}$, since, by assumption, $t$ is not a member of:
\[ f\fens{\LNPlat[\geqslant n]{\cF_U}{X}} \cap \Abhy{X} \]
a fortiori, it is not a member of:
\[ f\fens{\LNPlat[\geqslant n+1]{\cF_U}{X}} \cap \Abhy{X} \]
It then follows from \ref{THEO - Aplatissement fort cdt6*} that there exists a predecessor $v$ of $s$ in $V \times_Z S$, and thus in $\Lambda$.
\end{itemize}
Let $u$ be the image of $v$ in $Z$. Let’s show that $u$ is not a member of:
\[ \omega\fens{\Lambda} \cap \Abhy{Z} \]
We will use a proof by contradiction, assuming that $u$ is a limit point of $\omega\fens{\Lambda} \cap \Abhy{Z}$. Since $S$ is a divisor of $Z$, it does not intersect $\Abhy{Z}$, and therefore $u$ is a limit point of:
\[ \omega\fens{\LNPlat[\geqslant n]{\cF_V \oslash_Z S}{Z}} \cap \Abhy{Z \priv S} \]
Let $\zeta$ be a point of :
\[ \omega\fens{\LNPlat[\geqslant n]{\cF_V \oslash_Z S}{Z}} \cap \Abhy{Z \priv S} \]
By definition, the point $\zeta$ has a preimage $\eta$ via $\omega$ lying on $\LNPlat[\geqslant n]{\cF_V \oslash_Z S}{Z}$, and since $\zeta$ does not belong to $S$, the point $\eta$ in fact belongs to $\LNPlat[\geqslant n]{\cF_V}{Z}$. Now, since the morphism from $Z$ to $X$ is flat outside of $S$, this implies that the image of $\eta$ on $Y$ belongs to $\LNPlat[\geqslant n]{\cF}{X}$. Consequently, the image $\xi$ of $\zeta$ on $X$ belongs to $f\fens{\LNPlat[\geqslant n]{\cF}{X}}$. Furthermore, the point $\zeta$ belongs to $\Abhy{Z \priv S}$, which implies that $\xi$ belongs to $\Abhy{X}$ by property \ref{PROP - point ab par enjoliveur}. It follows from the above that $t$ belongs to the adherence of $f\fens{\LNPlat[\geqslant n]{\cF}{X}} \cap \Abhy{X}$, which is absurd. Therefore, the point $u$ does not belong to the adherence of:
\[ \omega\fens{\LNPlat[\geqslant n]{\cF_V \oslash_Z S}{Z}} \cap \Abhy{X} \]
The assertion \ref{THEO - Aplatissement fort cdtF} then ensures that $v$ belongs to the image of $V_1 \times_W \Sigma$ in $V$, and $s$ therefore belongs to the image of $V_1 \times_W \Sigma$ in $Y$.

\item \textit{Verification of the seventh assertion of the theorem:} \\
To demonstrate this, we proceed in two steps.
\begin{itemize}
\item \textit{Verification of the assertion} \ref{THEO - Aplatissement fort cdt7*} \textit{for} $i$ \textit{equal to} $N - (n + 1)$\\
Let $s$ be a point in $\Sigma \priv (S \cap W)$. There exists a point $u$ in:
\[ \LNPlat[\geqslant n]{\cF_{V \times_Z W}}{W} \]
such that $d_K(u) = n$ by \ref{THEO - Aplatissement fort cdtG}. Now, since the morphism from $W$ to $X$ is an isomorphism outside of $S \cap W$, the open set $U \times_X W \priv U \times_X (S \cap W)$ of $U \times_X W$ is an isomorphism onto its image in $U$, and $u$ belongs to $U \times_X W \priv U \times_X (S \cap W)$. Let $v$ denote the image of $u$ in $U$, we have:
\[ d_K(u) = d_K(v) \leqslant \dimAn{K}{\SuppFa{\cF_U}} \]
Let $t$ be the image of $s$ in $X$. It follows that $t$ belongs to $f\fens{\SuppFa{\cF_U}}$ and that:
\[ d_K(s) = d_K(t) \leqslant \dimAn{K}{\SuppFa{\cF_U}} - (N - (N - n - 1) - 1)) \]
\item \textit{Verification of the seventh assertion of the theorem:} \\
It follows from all of the above that there exists a local enjoliveur $(Z,S)$ of $(X,x)$ of class less than $N - n$ given by a factorization:
\[ \varphi : \xymatrix{ Z = Z_{2(N-n)} \ar[r] & Z_{2(N-n)-1} \ar[r] & \ldots \ar[r] & Z_0 = X } \]
and, for every $i$ between $0$ and $2(N-n)$, a closed $K$-analytic subspace $S_i$ of $Z_i$ such that:
\[ S_{2(N-n)} = S \quad\et\quad S_0 = \emptyset \]
and such that, for every $i$ between $0$ and $N-(n+1)$, the image of $S_{2i}$ on $X$ is contained in $f\fens{\SuppFa{\cF_U}}$ and, if $\cF_U$ is nonzero, this image has dimension less than:
\[ \dimAn{K}{\SuppFa{\cF_U}} - (N - (i - 1)) \]
Let $s$ be a point in $S$, since $S_0$ is empty and since, for every $i$ between $1$ and $N - n$, the set $S_{2i}$ is the inverse image of $S_{2i-1}$ under the action of $Z_{2i}$, there exists an $i$ between $0$ and $N - (n + 1)$ such that the image of $s$ in $Z_{2i+1}$ belongs to $S_{2i+1}$ but not to the inverse image of $S_{2i}$. This image is therefore located on:
\[ S_{2i+1} \priv (Z_{2i+1} \times_{Z_{2i}} S_{2i+1} ) \]
It follows from the above that the image $t$ of $s$ on $X$ belongs to $f\fens{\SuppFa{\cF_U}}$ and that we have:
\[ d_K(t) \leqslant \dimAn{K}{\SuppFa{\cF_U}} - (N - i - 1)) \leqslant \dimAn{K}{\SuppFa{\cF_U}} - n \]
\end{itemize}

\end{enumerate}
\end{enumerate}
\end{enumerate}
\end{dem}

\newpage \section{A non-Archimedean local Chevalley theorem in the boundaryless case}

\subsection{Locally embellishable morphisms}

\begin{att}
In this subsection, $(X,x)$ denotes a good reduced strictly pointed $K$-analytic space, $(Y,y)$ a good reduced strict pointed $K$-analytic space of dimension $n$, $f$ a pointed morphism from $(Y,y)$ to $(X,x)$ such that $y$ lies in the interior of $f$, $d$ the minimum value of the relative dimension of $f$ and $\delta$ the maximum value of the relative dimension of $f$.
\end{att}

\begin{deff}[Embellishable morphism at a point]
Let $k$ be an integer, we say that $f$ is $k$-embellishable at $y$ if there exists a local enjoliveur $(Z,S)$ of $(X,x)$ of class less than $k+1$ and an analytic neighborhood $U$ of $y$ in $Y$ such that, denoting by $V$ the $K$-analytic space $U \times_X Z$, we have:
\begin{enumerate}[label=$(EMB_{\arabic*})$]
\item The morphism from $V \oslash_Z S$ to $U$ is surjective.
\item The morphism from $V \oslash_Z S$ to $Z$ factors as a flat and surjective morphism onto a reduced closed analytic subspace $F$ of $Z$ of pure dimension $n-d$.
\item The image of $S$ in $X$ is contained in that of $U$ and has a dimension strictly less than $n-d$.
\end{enumerate}
If $k$ equals $\delta$, we simply say that $f$ is embellishable at $y$. Furthermore, we say that $f$ is embellished by $(Z,S,U,F)$ at $y$ and that $(Z,S,U,F)$ is a local embellisher of $f$ at $y$.
\end{deff}

\begin{lem}\label{LEM - embellissement sur un ferme}
Suppose there exists a closed and reduced analytic subspace $X_0$ of $X$ such that $f$ has a factorization of the form:
\[ f : \xymatrix{ (Y,y) \ar[rr]^g && (X_0,x) \ar[rr] && (X,x) } \]
If $g$ is embellishable at $y$, then $f$ is embellishable at $y$.
\end{lem}
\begin{dem}
Assume there exists a local enjoliveur $(W,R)$ of $(X_0,x)$ of class less than $\delta+1$ and an analytic neighborhood $U$ of $y$ in $Y$ such that, letting $V$ denote the $K$-analytic space $U \times_X W$, we have:
\begin{enumerate}
\item The morphism from $V \oslash_W R$ to $U$ is surjective
\item The morphism from $V \oslash_W R$ to $W$ factors as a flat and surjective morphism onto a reduced closed analytic subspace of $W$ of pure dimension $n-d$. 
\item The image of $R$ in $X_0$ is contained in the image of $U$ and has a dimension strictly less than $n-d$. 
\end{enumerate}
By the proposition \ref{PROP - remonter des enjoliveurs locaux}, there exists a local enjoliveur $(Z,S)$ of $(X,x)$ of class less than $\delta+1$ such that $W$ is identified with $(X_0 \times_X Z) \oslash_ZS$, such that $R$ is equal to $W \times_Z S$, and such that the image of $S$ in $X$ is equal to that of $R$, then we have:
\[ V \oslash_W R = V \oslash_Z S \]
The morphism $f$ is then embellishable at $y$.
\end{dem}

\begin{lem}\label{LEM - embellissement après enjoliveur}
Let $c$ and $c'$ be two integers whose sum is $\delta-d$, $(W,R)$ a local enjoliveur of $(X,x)$ of class $c$, and $U$ an analytic neighborhood of $y$ in $Y$ such that:
\begin{enumerate}
\item The morphism from $(U \times_X W) \oslash_W R$ to $U$ is surjective.
\item The image of $R$ in $X$ is contained in the image of $U$.
\item The relative dimension of the morphism from $(U \times_X W) \oslash_W R$ to $W$ is bounded above by $d+c'$.
\item The morphism from $(U \times_X W) \oslash_W R$ to $W$ is $(d+c')-$embellishable\footnote{By Lemma \ref{LEM - dimension et enjoliveurs locaux}, the morphism from $(U \times_X W) \oslash_W R$ is purely of dimension $n$, and the generic dimension of the morphism from $(U \times_X W) \oslash_W R$ into $W$ is equal to $d$. As for $W$, it is reduced.} at every point of:
\[ ((U \times_X W) \oslash_W R)_y \]
\end{enumerate}
If the set of points at which the relative dimension of $f$ is $d$ is dense, then $f$ is embellishable.
\end{lem}
\begin{dem}
Even if it means restricting $U$ to the neighborhood of $y$, we can assume that $U$ is closed in $Y$.
\begin{enumerate}
\item Let $V$ denote the $K-$analytic space $U \times_X W$. To say that the morphism from $V \oslash_W R$ to $W$ is embellishable at every point of $(V \oslash_W R)_y$ means that, for every point $t$ in $V \oslash_W R$ whose image in $Y$ is $y$ and whose image in $W$ is denoted by $w_t$, there exists a local enjoliveur $(Z_t,S_t)$ of $(W,w_t)$ of class $d+c'+1$ and an analytic neighborhood $U_t$ of $t$ in $V \oslash_W R$, which we may assume to be closed, such that, denoting by $V_t$ the $K-$analytic space $U_t \times_W Z_t$, we have:
\begin{enumerate}
\item The morphism $V_t \oslash_{Z_t} S_t$ into $U_t$ is surjective.
\item The morphism from $V_t \oslash_{Z_t} S_t$ to $Z_t$ factors as a flat and surjective morphism onto a reduced closed analytic subspace of $Z_t$ of pure dimension $n-d$.
\item The image of $S_t$ on $W$ is contained in the image of $U_t$ and has a dimension strictly less than $n-d$.
\end{enumerate}
Let $\cI_t$ be the coherent sheaf of ideals on $Z_t$ corresponding to $S_t$ and $\cJ_t$ be the one corresponding to $Z_t \times_W R$. Let $E_t$ denote the blow-up of $Z_t$ along $\cI_t \cJ_t$, and $D_t$ the corresponding exceptional divisor. Since the morphism from $Z$ to $W$ defines an enjoliveur, the fiber product $Z_t \times_W R$ is a divisor of $Z$ and, consequently, the ideal $\cI_t \cJ_t$ is invertible, so that $E_t$ identifies with $Z_t$ and $D_t$ with the sum of the divisors $S_t$ and $Z_t \times_W R$. We will now use the notation $Z_t$ instead of $E_t$, but retain the notation $D_t$. By construction, the pair $(Z_t,D_t)$ is a local enjoliveur of $(X,x)$ of class bounded by $d+c+c'+1 = \delta+1$.
Let $s$ be a point in $D_t$, we denote its image on $W$ by $w_t^s$ and its image on $X$ by $x_t^s$.
\begin{enumerate}
\item If $w_t^s$ belongs to $R$, then $x_t^s$ belongs to the image of $U$, and $d_k(x_t^s)$ is strictly less than $n-d$.
\item Otherwise, the point $s$ belongs to $S_t$, in which case $w_t^s$ belongs to the image of $U_t$ and satisfies the upper bound:
\[ d_k(w_t^s) < n-d \]
so that $x_t^s$ belongs to the image of $Y$ and $d_k(x_t^s)$ is strictly less than $n-d$.
\end{enumerate}
The analytic domain $V_t$ of $(V \oslash_W R) \times_W Z_t$ is the trace of a strictly compact analytic domain $\Omega_t$ of $V \times_W Z_t$. Now, an irreducible component of $\Omega_t$ is contained in the inverse image of $D_t$ if and only if it is contained in the inverse image of $Z_t \times_W R$ or of $S_t$, so $\Omega_t \oslash_Z D_t$ coincides with $V_t \oslash_{Z_t} S$.
\item The morphism $\varphi$ from $W$ to $X$, being a composition of domanial morphisms and blow-ups, is a topologically proper morphism. It follows that the morphism $\varphi_{V \oslash_W R}$ from $V \oslash_W R$ to $U$ is topologically proper by change of basis. The fiber $(V \oslash_W R)_y$ is then compact, and thus there exists a finite set of points $\Delta$ in $(V \oslash_W R)_y$ such that the interiors of the sets in the family $(U_t)_{t \in \Delta}$ form a covering of $(V \oslash_W R)_y$. Let $T$ denote the union of $(U_t)_{t \in \Delta}$, the set $T$ is a neighborhood of $(V \oslash_W R)_y$ in $V \oslash_W R$, and since $\varphi_V$ is topologically proper, the image of $T$ under $\varphi_V$ in $U$ is a neighborhood of $y$ in $Y$, and, without loss of generality, by restricting $U$ to the neighborhood of $y$, we may assume that $V \oslash_W R$ is equal to $T$.
\end{enumerate}
Let $U_1$ denote the disjoint sum of $(U_t)_{t \in \Delta}$, $V_1$ the disjoint sum of $(V_t)_{t \in \Delta}$, $\Omega_1$ the disjoint sum of $(\Omega_t)_{t \in \Delta}$, $Z_1$ the disjoint sum of $(Z_t)_{t \in \Delta}$, and $D$ the disjoint sum of $(D_t)_{t \in \Delta}$. Since the morphism from $V \oslash_W R$ to $U$ is surjective, it follows that the morphism from $\Omega \oslash_Z D$ to $U$ is surjective. Furthermore, the morphism from $\Omega \oslash_Z D$ to $Z$ factors through a flat and surjective morphism onto a reduced closed analytic subspace $F$ of $Z$ of pure dimension $n-d$. The local enjoliveur $(Z,D)$ of $(X,x)$ and the neighborhood $U$ of $y$ in $Y$ then define a local embellisher $(Z,D,U,F)$ of $f$ at $y$.
\end{dem}

\begin{theo}\label{THEO - embellissement de morphismes}
If the set of points at which the relative dimension of $f$ is $d$ is dense, then $f$ is embellishable at $y$ and, for every local embellisher $(Z,S,U,F)$ of $f$ at $y$, we have:
\begin{enumerate}
\item The space $(U \times_X Z) \oslash_Z S$ is reduced and purely $n$-dimensional.
\item The morphism $(U \times_X Z) \oslash_Z S$ into $U$ is generically an immersion of an analytic domain.
\item The morphism from $(U \times_X Z) \oslash_Z S$ to $F$ is purely of relative dimension $d$.
\item The set $S \cap F$ has an empty interior in $F$.
\end{enumerate}
\end{theo}
\begin{dem}\sautlistebis
\begin{enumerate}
\item \textit{Let's show that $f$ is embellishable at $y$} :\\
The proof will begin by reducing the problem to the relatively equidimensional case, and then proceed to the generically flat case.
\begin{enumerate}
\item \textit{Reduction when $\delta$ equals $d$:}\\
By applying Theorem \ref{THEO - Applatissement fort} to the coherent $\cOf{Y}$-module $\cOf{Y}$ and to the integers $\delta+1$ and $d+1$, we obtain the existence of a local enjoliveur $(Z,S)$ of $(X,x)$ of class less than $\delta - d$ and an analytic neighborhood $U$ of $y$ in $Y$ such that, letting $V$ denote the $K$-analytic space $U \times_X Z$, we have:
\begin{itemize} 
\item The morphism from $V \oslash_Z S$ to $Z$ is flat in dimensions greater than $d+1$. 
\item The image of $V \oslash_Z S$ in $U$ contains $U \priv \LNPlatS[\geqslant d+1]{U}{X}$.
\item The image of $S$ on $X$ is contained in the image of $U$ and has dimension less than $n-(d+1)$.
\end{itemize}
Since the morphism from $V\oslash_Z S$ to $Z$ is generically of dimension $d$ by Lemma \ref{LEM - dimension et enjoliveurs locaux}, the fact that this morphism is flat in dimensions greater than $d+1$ simply means, by (4.2) of \cite{Duc23}, that it is purely of relative dimension $d$. It also follows that $\LNPlat[\geqslant d+1]{Y}{X}$ is the locus of points in $Y$ at which the relative dimension of $f$ is $d+1$. The image of $V\oslash_Z S$ in $U$ therefore contains all fibers of $\rest{f}{U}$ that are purely of dimension $d$, and in particular, by (4.2) from \cite{Duc23}, all points of $\Abhy{U}$. This image is therefore dense in $U$, and since it is compact, it is the entire set $U$. It follows from Lemma \ref{LEM - embellissement après enjoliveur}, applied to $\delta-d$ and $0$, that it suffices to prove the theorem for the morphism from $V\oslash_Z S$ to $Z$. We have thus reduced the problem to the case where $\delta$ equals $d$, that is, once again to the case where the morphism $\rest{f}{U}$ is purely of relative dimension $d$.
\item \textit{Reduction to the generically flat case:}\\
Suppose that $f$ is purely of relative dimension $d$.
\begin{enumerate}
\item It follows, by virtue of the property \ref{PROP - factorisation fini et lisse connexe}, that there exists a strictly affinoid rational neighborhood $V$ of $y$ in $Y$, a strictly affinoid neighborhood $W$ of $x$ in $X$, and a $d$-polyray $r$ with coefficients in $\abs{K^\times}\!{}^{\eQ}$ such that $\rest{f}{V}$ factors through a finite morphism from $V$ to $T = \eD^d_W(0,r)$. Since the morphism from $V$ to $T$ is finite, its image $F$ is a Zariski closed set in $T$, and since $\rest{f}{V}$ is purely of relative dimension $d$ and the fibers of the projection of $T$ onto $W$ are irreducible and of dimension $d$, the closed set $F$ is a union of fibers of $\rest{f}{V}$. Let $\omega$ be the morphism from $T$ to $W$, the existence of an $\omega$-coefficient ideal for $F$, guaranteed by Theorem (3.9) of \cite{Duc20}, then ensures that the image of $F$ in $W$ is a Zariski closed set $\Lambda$ of $W$. The morphism from $F$ to $\Lambda$, with the source and target endowed with their reduced structures, is quasi-smooth and purely of relative dimension $d$. Furthermore, the space $F$ is purely $n$-dimensional because it is the image of $V$ under a finite morphism and because $V$ itself is purely $n$-dimensional. It follows that $\Lambda$ is purely $n-d$-dimensional.
\item It suffices to prove the desired result for the morphism from $V$ to $W$ and, consequently, by Lemma \ref{LEM - embellissement sur un ferme}, for the morphism from $V$ to $\Lambda$. Indeed, suppose there exists a local enjoliveur $(Z,S)$ of $(W,x)$ of class less than $d+1$ and an analytic neighborhood $U$ of $y$ in $V$ such that, denoting $V_1$ the $K$-analytic space $U \times_W Z$, we have:
\begin{itemize}
\item The morphism from $V_1 \oslash_Z S$ to $U$ is surjective.
\item The morphism from $V_1 \oslash_Z S$ to $Z$ factors as a flat and surjective morphism onto a reduced closed analytic subspace $R$ of $Z$ of pure dimension $n-d$.
\item The image of $S$ in $W$ is contained in the image of $U$ and has a dimension strictly less than $n-d$. 
\end{itemize}
Let us choose a representation of the local enjoliveur $(Z,S)$ of $(W,x)$:
\[ \xymatrix{ Z = Z_{2(d+1)} \ar[rr] && Z_{2d+1} \ar[rr] && \ldots \ar[rr]^u && Z_0 = W } \]
where, for all $i$ between $0$ and $2n$, a closed $K$-analytic subspace $S_i$ of $Z_i$ is given such that $S_{2(d+1)} = S$ and $S_0 = \emptyset$.
By replacing $Z_0$ with $X$ and $u$ with the composition of $u$ and the inclusion of $W$ into $X$, we obtain a local enjoliveur $(Z,S)$ of $(X,x)$ of length less than $d+1$. It is then immediately clear that $(Z,S,U,R)$ is a local embellisher of $f$ at $y$.
\end{enumerate}
It is therefore sufficient to prove the theorem for the morphism from $V$ to $\Lambda$, which reduces the problem to the case\footnote{The morphism from $V$ to $\Lambda$ is also relatively isodimensional, but this is no longer relevant at this stage of the proof, so we disregard this assumption.} where $X$ is purely of dimension $n-d$.
\item \textit{Proof for the generic flat case:}\\
Suppose that $X$ is purely $n - d$-dimensional. In this case, the image of $\Abhy{Y}$ is contained in $\Abhy{X}$, and the Zariski open set $\LPlat{Y}{X}$ is dense in $Y$ and contains the entire fiber of every point in $\Abhy{Y}$ by (4.2) in \cite {Duc23}. Applying Theorem \ref{THEO - Applatissement fort} to the coherent $\cOf{Y}$-module $\cOf{Y}$ and to the integers $\delta+1$ and $0$, there exists a local enjoliveur $(Z,S)$ of $(X,x)$ of class less than $\delta+1$ and an analytic neighborhood $U$ of $y$ in $Y$ such that, letting $V$ denote the $K$-analytic space $U \times_X Z$, we have:
\begin{enumerate}
\item The morphism from $V \oslash_Z S$ to $Z$ is flat.
\item The image of $V \oslash_Z S$ in $U$ contains:
\[ U \priv \LNPlatS[\geqslant 0]{U}{X} \]
This is equivalent to saying that the image of $V \oslash_Z S$ in $U$ contains the union of the fibers of $\rest{f}{U}$ that are entirely contained in $\LPlat{U}{X}$.
\item The image of $V \times_Z S$ in $U$ is contained in $\LNPlat[\geqslant 0]{U}{X}$, which means that it does not intersect any fiber of $\rest{f}{U}$ that is entirely contained in $\LPlat{U}{X}$.
\item The image of $S$ in $X$ is contained in the image of $U$ in $X$ and has dimension less than $n$.
\end{enumerate}
Since the morphism from $V \oslash_Z S$ to $Z$ is flat, its image $Z_0$ on $Z$ is a strictly compact analytic domain of $Z$. We can then replace $Z$ with $Z_0$ and $S$ with $S \cap Z_0$ without altering the above properties, which allows us to assume that the morphism from $V \oslash_Z S$ to $Z$ is surjective. The image of $V \oslash_Z S$ in $U$ is a compact subset of $U$ containing the union of the fibers of $\rest{f}{U}$ that are entirely included in $\LPlat{U}{X}$. In particular, it contains $\Abhy{U}$, which is dense in $U$. Consequently, the morphism from $V \oslash_Z S$ to $U$ is surjective. Let $s$ be a point in $S$ and $t$ its image in $X$, it remains to verify that we have:
\[ d_k(t) \leqslant n-d-1 \]
Since the morphism from $V \oslash_Z S$ to $U$ is surjective, the point $s$ has a predecessor $v$ in $V$, whose image $u$ in $U$ is a predecessor of $t$. However, since $u$ is the image of an element of $V \times_X S$, its fiber is not entirely contained in $\LPlat{U}{X}$, and this rules out the possibility that $t$ belongs to $\Abhy{X}$ by (4.2) of \cite{Duc23}. Consequently, we have:
\[ d_k(t) <  n - d \]
\end{enumerate}

\item \textit{Let's prove the rest of the assertion:}\\
It follows from the definition of local embellishments that the morphism from $V \oslash_Z S$ to $U$ is generically an isomorphism and that $V \oslash_S Z$ is reduced and purely $n$-dimensional. The latter fact follows from the fact that $Y$ itself is purely $n$-dimensional. It follows from the triviality of the morphism from $V \oslash_Z F$ to $F$ that its fibers are purely $n-d$-dimensional. Finally, let $s$ be a point in $S \cap F$, the fiber of $V \oslash_Z S$ at $s$ is nonempty by the surjectivity of $V \oslash_Z S$ onto $F$, and it is purely $d$-dimensional. It therefore has a point $t$ such that $d_{\CrpResC{s}}(t)$ equals $d$. But by the definition of $V \oslash_Z S$, the reciprocal image of $S$ under this map has an empty interior. Consequently, the integer $n$ strictly greater than $d_K(t)$ is equal to $d_k(s) + d$, and $d_k(s)$ is strictly less than $n - d$. It follows that $S \cap F$ has an empty interior in $F$.
\end{enumerate}
\end{dem}

\subsection{A non-Archimedean local Chevalley theorem in the boundaryless case}

\begin{nott}
Let $X$ be a topological space and $C$ a locally closed subset of $X$; we denote $\partial C$ as the complement of $C$ in its closure in $X$.
\end{nott}

\begin{deff}[Cell]
Let $X$ be a $K$-analytic space, $C$ a locally closed subset of $X$, and $n$ an integer. We will define recursively in terms of $n$ the property that $C$
is an $n$-cell.
\begin{enumerate}
\item $C$ is said to be a $0-$cell if $C$ is a nonempty Zariski closed set of $X$ of dimension zero. 
\item If $n$ is nonzero, $C$ is said to be an $n$-cell if: 
\begin{enumerate}
\item There exists a finite family of closed analytic sets $(V_i)_{i \in I}$ in $X \priv \partial C$ and a nonempty Zariski set $Z$ consisting of the disjoint sum of purely $n$-dimensional sets $(V_i)_{i \in I}$ whose image in $X$ is $C$.
\item There exists a finite family of integers $(n_j)_{j \in J}$ strictly less than $n$ and a family $(C_j)_{j \in J}$ of subsets of $X$ such that, for every $j$ in $J$, the set $C_j$ is an $n_j$-cell and such that $\partial C$ is the union of $(\adh{C_j})_{j \in J}$.
\end{enumerate}
\end{enumerate}
\end{deff}

\begin{theo}\label{Chevalley non archimedien}
Let $(X,x)$ be a compact pointed $K$-analytic space, $(Y,y)$ a compact equidimensional $K$-analytic space of dimension $n$ and $f$ a pointed morphism from $(Y,y)$ to $(X,x)$, if there exists an integer $d$ such that $f$ is generically of dimension $d$, then there exists an analytic neighborhood $U$ of $y$ in $Y$, an analytic neighborhood $V$ of $x$ in $X$, and an $(n-d)$-cell $C$ of $V$ such that:
\begin{enumerate}
\item The image of $U$ under $f$ is equal to $\adh{C}$.
\item The set $f\fens{\Abhy{U}}$ is contained in $C$.
\end{enumerate}
\end{theo}
\begin{dem}\sautlistebis
\begin{enumerate}
\item Let us prove the first assertion by strong recurrence in dimension $n-d$. We therefore assume the result to be true in dimensions strictly less than $n-d$. Furthermore, since the statement of the theorem is invariant under nilpotents, we may assume that $X$ and $Y$ are reduced. Theorem \ref{THEO - embellissement de morphismes} then guarantees that there exists a local embellisher $(Z,S,U,F)$ of $f$ at $y$ that satisfies the following properties:
\begin{enumerate}
\item The morphism from $V \oslash_Z S$ to $U$ is surjective.
\item The morphism from $V \oslash_Z S$ to $Z$ factors as a flat and surjective morphism onto a reduced closed analytic subspace $F$ of $Z$ of pure dimension $n-d$.
\item The image of $S$ in $X$ is contained in that of $U$ and has a dimension strictly less than $n-d$.\label{CDT c chev}
\item The set $S \cap F$ has an empty interior in $F$.
\end{enumerate}
Let $\varphi$ denote the morphism from $Z$ to $X$. Since the morphism from $V \oslash_Z S$ to $U$ is surjective, the set $f\fens{U}$ is equal to the image of $V \oslash_Z S$ in $X$, that is, to $\varphi\fens{F}$. Let $\Sigma$ denote the set $\reci{\varphi}\fens{\varphi\fens{S}}$; this is a closed set in $Z$, and we have:
\[ f\fens{U} = \varphi\fens{F} = \varphi\fens{F \priv \Sigma} \sqcup \varphi\fens{F \cap \Sigma} = \varphi\fens{F \priv \Sigma} \sqcup \varphi\fens{S} \]
where the last equality follows from the fact that $\varphi\fens{S}$ is contained in $f\fens{U}$. Let $z$ be a point in $\Abhy{F}$, it is not on $S$, so $\varphi$ is quasi-étale at $z$. In particular, $d_K(\varphi(z))$ equals $d_K(z)$, that is, $n-d$, which ensures that $\varphi(z)$ does not belong to $\varphi\fens{S}$. Consequently, $\Sigma$ does not intersect $\Abhy{F}$. Let $\psi$ denote the morphism from $Z \priv \Sigma$ to $X \priv \varphi\fens{S}$ induced by $\varphi$ via a change of basis along the open immersion $X \priv \varphi\fens{S}$ into $X$. The morphism $\psi$ is compact and is a composition of domanial morphisms and blow-ups. Let $C$ be the image of $F \priv\Sigma$ under $\varphi$. Let $z$ be a point in $\Abhy{F}$. We saw earlier that $z$ does not belong to $\Sigma$, and thus $\varphi(z)$ belongs to $C$. Since $\Abhy{F}$ is dense in $F$, this implies that $C$ is dense in $\varphi\fens{F}$. It therefore suffices to show that $C$ is an $(n-d)$-cell to conclude the result. Given its definition, the only point that remains to be verified is that $\partial C$ can be written as a finite union of local cell adhesions of appropriate dimensions. Since $C$ is dense in $\varphi\fens{F}$, we have $\partial C$ equals $\varphi\fens{S}$. For any pair $e$ of integers $e_1$ and $e_2$ such that $e_2$ is less than $e_1$, let $S_e$ denote the union of the irreducible components of $S$ of dimension $e_1$ whose relative generic dimension over $X$ is equal to $e_2$. If $e$ is such that $S_e$ is nonempty, then, by \ref{CDT c chev}, that the dimension of the image of $S_e$ under $\varphi$ is strictly less than $n-d$, and our recurrence hypothesis then ensures that $\varphi\fens{S_e}$ is the adherence of a cell $D_e$ of dimension equal to the dimension of the image of $S_e$ under $\varphi$. Let $E$ denote the finite set of indices $e$ such that $S_e$ is nonempty; then $S$ is the union of $(S_e)_{e \in E}$, and $\partial C$ is the union of the adheres of $(D_e)_{e \in E}$.

\item First, note that for any $y$ in $\Abhy{U}$, the integer $d_K(f(y))$ is equal to $n-d$. Therefore, if $y$ is a point in $Y$, $f(y)$ cannot lie in $\partial C$, so $f(y)$ belongs to $C$, which implies the second statement.
\end{enumerate}
\end{dem}

\newpage
\bibliographystyle{perso}
\bibliography{bibli}
\end{document}